УДК 519.217.2

*Посвящается Виталию Дмитриевичу Арнольду (4.10.1968 – 04.01.2017)*

# Вокруг степенного закона распределения компонент вектора PageRank[1]

*Александр Владимирович Гасников (МФТИ, ИППИ РАН), Максим Евгеньевич Жуковский (МФТИ, Яндекс),*

*Сергей Вячеславович Ким (2007 школа, г. Москва), Федор Андреевич Носков (Лицей № 2, г. Тула),*

*Степан Сергеевич Плаунов (2007 школа, г. Москва), Даниил Андреевич Смирнов (2007 школа, г. Москва)*

**Аннотация**

В данной статье с помощью большого числа компьютерных экспериментов обнаруживается степенной закон распределения компонент вектора PageRank, посчитанного для web-графа, сгенерированного по модели Бакли–Остгуса. Подробно обсуждается выбор алгоритмов, которые использовались в экспериментах по расчету вектора PageRank. Рассматривается приближенная к реальности модель ранжирования web-страниц, предложенная совместно с сотрудниками компании Яндекс. Обсуждаются вычислительные аспекты предложенной модели в контексте ранее описанных численных способов поиска вектора PageRank.

**Ключевые слова:** марковская цепь, эргодическая теорема, мультиномиальное распределение, концентрация меры, оценка максимального правдоподобия, Google problem, градиентный спуск, автоматическое дифференцирование, степенной закон.

## 1. Введение

В данной статье обсуждается несколько фундаментальных вопросов, связанных с: эргодической теоремой для марковских процессов, методами Монте-Карло, явлением концентрации меры, понятием равновесия макросистемы, основной теоремой математической статистики (теорема Фишера) о свойствах оценки максимального правдоподобия, ролью степенных законов распределения, машинным обучением, невыпуклой оптимизацией, автоматическим дифференцированием и рядом смежных вопросов. Все темы/вопросы обсуждаются на одном примере – Google problem. Умение эффективно и многократно решать эту задачу, например,

---

позволяет поисковым системам ранжировать web-страницы в ответ на пользовательский запрос [42, 51, 53, 62, 72, 87]. Схожие постановки задач возникают в моделях поиска консенсуса [1]. Собранные в статье математические теоремы (конструкции), безусловно, входят в "золотой фонд" современной математики. Тем интереснее, как нам кажется, будет узнать о том, как отмеченные конструкции сосуществуют. Близкие сюжеты собраны также в пишущейся сейчас книге [50].

## 2. Google problem и эргодическая теорема

В 1998 г. Ларри Пейджем и Сергеем Брином был предложен специальный способ ранжирования web-страниц, который и лег в основу поисковой системы Google [72, 53]. Далее мы опишем этот способ, предварив описание необходимыми сведениями из теории дискретных марковских цепей [10, 29, 47, 58, 61, 70, 74–76, 78, 93].

Рассмотрим ориентированный взвешенный граф (см. рис. 1). Граф имеет $n$ вершин. Ребро, выходящее из вершины $i$ в вершину $j$, имеет вес $p_{ij} \geq 0$. Если из вершины $i$ в вершину $j$ ребра нет, то полагаем $p_{ij} = 0$. Число $p_{ij}$ интерпретируется как вероятность перейти из вершины $i$ в вершину $j$. Поскольку распределение вероятностей должно быть нормировано на единицу, то для любой вершины $i$ имеет место равенство $\sum_{j=1}^{n} p_{ij} = 1$. Набор чисел $\{p_{ij}\}_{i,j=1,1}^{n,n}$ удобно будет записать в виде матрицы $P = \|p_{ij}\|_{i,j=1,1}^{n,n}$ (квадратной таблицы на пересечении $i$ строки и $j$ столбца которой стоит $p_{ij}$).

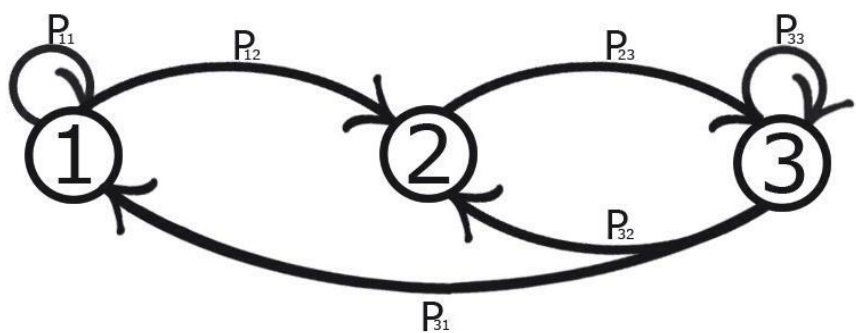

Рис. 1 Пример web-графа

На граф можно посмотреть как на город, вершинами которого являются различные районы города, а ребра – дороги (вообще говоря, с односторонним движением). Предположим, что в городе имеется "Красная площадь" – такой район, в который можно попасть по дорогам из любого другого. Оказывается, в этом предположении верен следующий результат (*эргодическая теорема для марковских процессов*): *если пустить блуждать человека по городу в течение длительного времени так, что человек будет случайно перемещаться из района в* *район согласно весам ребер графа, то доли* $\{v_k\}_{k=1}^{n}$ *времени, которые человек провел в разных районах будут удовлетворять следующему уравнению* $\sum_{i=1}^{n} v_i p_{ij} = v_j$ *или в векторном виде* $v^T P = v^T$, *имеющему единственное решение* (в классе $\sum_{i=1}^{n} v_i = 1$). В связи с последним представлением говорят, что[2] $v = (v_1, ..., v_n)^T$ является левым собственным вектором (стохастической по строкам) матрицы $P$. В экономике (теории неотрицательных матриц [36]) этот вектор также называют *вектором Фробениуса–Перрона*. Собственно, предположение о "Красной площади" и обеспечивает единственность $v$ (инвариантного / стационарного распределения). Если это предположение не верно, то вектор $v$ существует (то есть на больших временах возникают предельные пропорции), но он не единственен. И ответ (выбор $v$) уже будет зависеть от того, откуда стартовал человек. В общем случае город распадается на отдельные несвязанные между собой (существенные) кластеры[3] из районов, связанных между собой внутри кластера, и отдельного кластера несущественных районов. Этот кластер определяется тем, что из любого его района можно уйти в один из существенных кластеров, при этом попасть (обратно) в несущественный кластер из существенных кластеров невозможно (нет дорог). Из такого несущественного кластера человек в конечном итоге "свалится" в один из существенных кластеров, откуда уже не выйдет. Рис 2 поясняет общий случай.

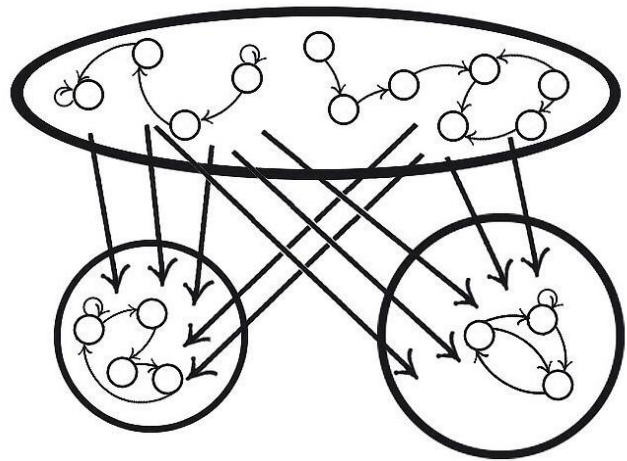

Рис. 2 Граф однородной дискретной марковской цепи (общий случай)

Чтобы понять, откуда получается приведенная выше эргодическая теорема (в вольной трактовке), получим с помощью формулы полной вероятности уравнение Колмогорова–Чэмпена (основное уравнение, описывающее *дискретные однородные марковские цепи*). Обо-значим через $p_i(t)$ – вероятность того, что человек находится в момент времени $t$ в районе $i$. Тогда по формуле полной вероятности для любого $j = 1,\ldots,n$

$$p_j(t+1) = \sum_{i=1}^{n} \underbrace{P\begin{pmatrix}\text{человек в момент времени } t \\ \text{находился в районе номер } i\end{pmatrix}}_{p_i(t)} \cdot \underbrace{P\begin{pmatrix}\text{человек перешел в район } j \\ \text{при условии, что был в районе } i\end{pmatrix}}_{p_{ij}}.$$

---

[2] Здесь и далее, если $A$ – матрица из $m$ строк и $n$ столбцов с элементами $A_{ij}$ – то, что стоит на пересечении $i$ строки и $j$ столбца, то по определению $A^T$ – матрица из $n$ строк и $m$ столбцов с элементами $A_{ji}$ – то, что стоит на пересечении $i$ строки и $j$ столбца. В частности, для рассматриваемого здесь случая $m = 1$

[3] То есть из каждого такого кластера в любой другой кластер дорог нет.

Или в матричном виде

$$p^T(t+1) = p^T(t)P. \qquad (1)$$

Последнюю формулу и называют **уравнением Колмогорова–Чэмпена**.

Предположим теперь, что существует предел $\lim_{t\to\infty} p(t) = v$. Какому уравнению должен удовлетворять вектор $v$? Переходя к пределу в обеих частях равенства (1) и учитывая, что $\lim_{t\to\infty}(p^T(t)P) = \lim_{t\to\infty}(p^T(t))P$, получим уже известное нам соотношение

$$v^T P = v^T. \qquad (2)$$

Именно из такого соотношения и было предложено искать вектор ранжирования web-страниц (**PageRank вектор**) в модели Брина–Пейджа. Отличие этой модели от описанной выше только в интерпретации: вершины – web-страницы, ребра – линки (гиперссылки).

Заметим, что в общем случае эргодическая теорема – это похожий по формулировке факт для общих динамических систем (у которых динамика в отличие от (1) существенно нелинейная, равно как и получаемый при этом аналог соотношения (2)). Грубо, эту теорему можно сформулировать так: *подобно (2) ищется инвариантная мера, т.е. такое распределение, которое не изменяется при действии оператора перехода динамической системы (в нашем случае оператор описан (1)), если это распределение единственно, то доли времени, которая динамическая система проведет в том или ином состоянии, описываются этой инвариантной мерой*. С примерами применения этой теоремы любознательный читатель мог встречаться в теории чисел: теорема Вейля [13], теорема Гаусса–Кузьмина [3], а также в курсе термодинамики [28] (эргодическая гипотеза Лоренца). На самом деле связь между только что разобранной конструкцией (эргодической теорией дискретных однородных цепей Маркова) и эргодической теорией динамических систем не ограничивается отмеченной аналогией. Связь эта намного глубже. В основе этой связи лежит *конструкция Улама* [7] (см. также [44]). К сожалению, все это существенно выходит за рамки даже институтских курсов, поэтому здесь мы лишь ограничимся ссылкой на мультфильм, который на наш взгляд хорошо поясняет написанное [95].

Стоит оговорить один нюанс. Приведенные выше рассуждения справедливы в предположении существования предела $\lim_{t\to\infty} p(t) = v$. Казалось бы, что предположения о "Красной площади" будет достаточно и тут. Однако, как показывает простейший пример (рис. 3), в котором $n = 2$, $p_{11} = p_{22} = 0$, $p_{12} = p_{21} = 1$, хотя вектор $v = (1/2, 1/2)^T$ существует и единственен (других решений в классе распределения вероятностей у (2) нет), тем не менее, предел $\lim_{t\to\infty} p(t)$ не существует, поскольку с ростом $t$ будет происходить периодическое чередование нулей и единиц в каждой компоненте вектора $p(t)$.

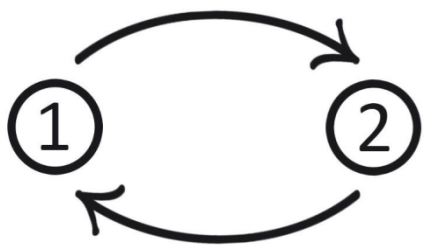

Рис. 3 Периодическая Марковская цепь (период равен 2)

Оказывается, что если выполняется условие "непериодичности", то предел, $\lim_{t\to\infty} p(t) = v$ действительно, существует. Более того, эти два условия (существования "Красной площади" и "непериодичности") являются не только достаточными, но и необходимыми для существования предела. Опишем, в чем заключается условие "непериодичности". Из "Красной площади" выходит много различных маршрутов, которые в конце снова приводят на "Красную площадь". Условие "непериодичности" означает, что наибольший общий делитель последовательности длин всевозможных маршрутов (начинающихся и заканчивающихся на "Красной площади") равен 1. Уточним, что длина маршрута равна числу ребер, вошедших в маршрут. В типичных web-графах оба отмеченных условия выполняются, поэтому в дальнейшем мы уже не будем делать соответствующие оговорки.

### 3. Стандартные численные подходы к решению Google problem

В предыдущем пункте мы не только описали Google problem, на самом деле, мы еще и описали, явно на это не указывая, один из основных численных методов поиска вектора PageRank [72]. Вспомним, что поиск этого вектора, сводится к решению задачи (2). С другой стороны, решение задачи (2) можно понимать как предел $v = \lim_{t\to\infty} p(t)$, где $p(t)$ рекуррентно рассчитывается согласно (1). Естественно, возникает желание использовать легко программируемую процедуру (1) для приближенного расчета вектора PageRank. Основной вопрос здесь: насколько быстро и точно имеет место сходимость в (1).

Назовем *спектральной щелью* матрицы $P$ наибольшее такое число $\alpha = 1 - |\beta| > 0$, где (вообще говоря, комплексное) число $\beta$ удовлетворяет условию: $|\beta| < 1$ и существует такой вектор $\eta$, что $\eta^T P = \beta \eta^T$ (хотя бы одно такое $\beta$ существует). Другими словами, $\alpha$ – расстояние между максимальным собственным значением матрицы $P$ (для стохастических матриц всегда равно 1) и следующим по величине модуля. Отсюда и название – спектральная щель (спектральный зазор – от англ. spectral gap). Оказывается, что имеет место следующий результат

$$\|p(t) - \nu\|_1 \overset{def}{=} \sum_{k=1}^{n} |p_k(t) - \nu_k| \le C \exp(-\alpha t / \tilde{C}). \quad (3)$$

Здесь и далее константы $C$, $\tilde{C}$ будут обозначать некоторые (каждый раз свои) универсальные константы, которые зависит от некоторых дополнительных деталей постановки, и которые, как правило, ограничены числом 10.

Заметим, что часто факт сходимость процесса (1) к решению (2) также называют эргодической теоремой, см., например, [29]. В прошлом пункте мы видели, что это, действительно, связанные вещи. Однако, на наш взгляд, сходимость процесса (1) – это скорее принцип сжимающих отображений (можно показать, что этот принцип, действительно, применим к (1), если под пространством понимать всевозможные лучи неотрицательного ортанта, а под метрикой на этом пространстве, в которой матрица $P$ "сжимает", понимать метрику Биркгофа–Гильберта [31]). К сожалению, детали здесь не совсем элементарны, поэтому мы ограничимся только следующей аналогией. Формула (3) отражает геометрическую скорость сходимости, характерную для принципа сжимающих отображений, а коэффициент $\alpha$ как раз и характеризует степень сжатия, осуществляемого матрицей $P$. Геометрически себе это можно представлять (правда, не очень строго) как сжатие с коэффициентом не меньшем $1 - \alpha$ к инвариантному направлению, задаваемому вектором $\nu$.

К сожалению, оценка (3) всего лишь переводит задачу оценивания скорости сходимости (1) в задачу оценивания спектральной щели $\alpha$. Однако последняя задача в ряде случаев довольно эффективно решается. В частности, к эффективным инструментам оценки $\alpha$ относятся изопериметрическое неравенство Чигера и неравенство Пуанкаре [29, 47, 61, 74, 93], которое также можно понимать как неравенство концентрации меры [93] (см. далее). Имеются и другие способы оценки $\alpha$ (см., например, [58, 75, 76, 78, 88]), однако все эти способы далеко выходят за рамки даже институтских курсов. Поэтому здесь мы ограничимся простыми, но важными в контексте рассматриваемых приложений, случаями.

А именно, предположим сначала, что для рассматриваемого web-графа существует такая web-страница, на которую есть ссылка из любой web-страницы в том числе из самой себя (усиленный аналог предположения о наличии "Красной площади" – поскольку теперь на красную площадь из любого района есть прямые дороги), более того предположим, что на каждой такой ссылке стоит вероятность не меньшая, чем $\gamma$. Для такого web-графа имеет место неравенство $\alpha \ge \gamma$ [83].

Предположим далее, что в модели блуждания по web-графу имеется "телепортация": с вероятностью $1-\delta$ человек действует так, как действовал бы раньше, а с вероятностью $\delta$ "забывает про все правила" случайно равновероятно выбирает среди $n$ вершин одну, в которую и переходит (возможно, что выберет и ту, в которой находился в момент выбора). Тогда, вводя квадратную матрицу $E$ размера $n$ на $n$, состоящую из одинаковых элементов $1/n$, перепишем уравнение (2) следующим образом

$$p^T(t+1) = p^T(t)((1-\delta)P + \delta E). \qquad (4)$$

В таком случае вектор PageRank необходимо будет искать из уравнения

$$v^T = v^T((1-\delta)P + \delta E). \qquad (5)$$

При $0 < \delta < 1$ уравнение (5) гарантированно (не зависимо от стохастической матрицы $P$) имеет единственное (в классе распределений вероятностей) решение. Более того, для (4), (5) имеет место оценка $\alpha \geq \delta$ [72]. На практике для подсчета PageRank обычно используют модель (5) с $\delta = 0.15$ [42, 72].

Поскольку матрица $P$ для реальных web-графов обычно сильно разреженная, то использовать формулу (4) в таком виде не рационально. Получается, что одна итерация (расчет по формуле (4)) будет стоить $3n^2$ арифметических операций (типа умножения чисел). Однако формулу (4) можно переписать следующим образом (раскрыли в (4) скобки, выполнив соответствующие умножения)

$$p^T(t+1) = (1-\delta)p^T(t)P + \delta \cdot (1/n, ..., 1/n). \qquad (5')$$

Последняя формула требует по порядку лишь $2sn$ арифметических операций, где $s$ – среднее число ненулевых элементов в строке матрицы $P$. Для реального интернета $n \approx 10^{10}$, а $s \ll 10^4$.

Только что был описан исторически самый первый алгоритм, использовавшийся для расчета вектора PageRank. Он получил название метода простой итерации (МПИ). Таким образом, МПИ эффективен, если $\alpha$ не очень близко к нулю. В частности, как в примере с телепортацией. Действительно, современный ноутбук в состоянии выполнять до $10^{10}$ арифметических операций в секунду (с учетом программной реализации и ряда других ограничений этот порядок на практике обычно уменьшается до $10^8$). Из формулы (3) видно, что сходимость, например, при $\alpha \geq 0.15$ очень быстрая. Нужная точность получается уже на нескольких десятках итераций. С учетом того, что одна итерация "стоит" по порядку $2sn \leq 10^{13}$ получаем программу, которая в состоянии за день найти для реального интернета вектор PageRank. К сожалению, на самом деле все не так просто. Проблема в памяти, в которую необходимо загружать матрицу $P$. Разумеется, $P$ необходимо загружать ни как матрицу из $n^2$ элементов, а в виде, так называемых, списков смежностей по строкам (для ряда других алгоритмов, еще и по столбцам). Однако это все равно не решает проблемы. Быстрая память компьютера – Кэш память (разных уровней). Кэш память (процессора) совсем маленькая, ее с огромным запасом ни на что такое не хватит. Более медленная – оперативная память. Тем не менее, если удалось бы выгрузить матрицу $P$ в такую память, то производительность программы по-прежнему соответствовала бы своим расчетным мощностям (см. выше). Обычно оперативной памяти в современном персональном компьютере (тем более ноутбуке) не более нескольких десятков Гигабайт ($\sim 10^{10}$ байт), чего, очевидно, не хватает. Следующая память – жесткий диск. Обращение программы к этой памяти, по-сути, останавливает нормальную работу программы. Как только кончаются ресурсы оперативной памяти, мы видим, что программа либо не работает совсем, либо начинает работать очень медленно. Отмеченную проблему можно решать, увеличивая оперативную память, либо используя распределенную память. Мы не будем здесь на этом останавливаться. Укажем лишь, что во избежание указанных проблем с памятью, все численные эксперименты в статье проводятся для размеров графов $n \leq 10^5$ (см. п. 7).

Но что делать, если $\alpha$ оказалось достаточно малым или мы не можем должным образом оценить снизу $\alpha$, чтобы гарантировать быструю сходимость МПИ? В таком случае полезным оказывается следующий простой результат [38] (Поляк–Тремба, 2012), не предполагающий, кстати, выполнения условия "непериодичности",

$$\left\| P^T \bar{p}_T - \bar{p}_T \right\|_1 \leq \frac{C}{T}, \quad \bar{p}_T = \frac{1}{T} \sum_{t=1}^{T} p(t). \tag{6}$$

Эта оценка уже никак не зависит от $\alpha$. Вектор $\bar{p}_T$ считается по порядку величины за то же время что и $p(T)$ (просто параллельно в процедуре (1) введем суммирование получающихся векторов). К сожалению, для наших целей: исследование степенного закона убывания компонент вектора PageRank метод Поляка–Трембы не подходит. Причина связана с видом оценки (6). Имеется принципиальная разница с оценкой (3) в том, что в оценке (3) мы можем гарантировать близость найденного вектора $p(t)$ к вектору PageRank $v$, обеспечив малость $\left\| p(t) - v \right\|_1$. Что касается соотношения (6), то в типичных случаях из $\left\| P^T \bar{p}_T - \bar{p}_T \right\| \approx \varepsilon$ можно при больших $n$ лишь получить, что $\left\| \bar{p}_T - v \right\| \approx C\varepsilon / \alpha$, т.е. опять возникает "нехорошая" зависимость в оценке от $\alpha$. Для симметричной матрицы $P$ и 2-нормы этот результат был получен Красносельским–Крейном в 1952 г. [30].

К сожалению, и другие способы поиска вектора PageRank, которые на первый взгляд не используют в своих оценках $\alpha$, на деле оказываются методами, выдающими такой вектор $\tilde{p}_T$, что $\left\| P^T \tilde{p}_T - \tilde{p}_T \right\| \approx \varepsilon$ в некоторой норме (обычно это 1-норма, 2-норма и бесконечная норма – в следующем пункте мы поясним, что имеется в виду под этими нормами), что сводит ситуацию к рассмотренной выше. Тем не менее, хочется отметить большой прогресс, достигнутый за последнее время (во многом благодаря разработкам Б.Т Поляка, А.С. Немировского, Ю.Е. Нестерова, см., например, [16] и цитированную там литературу), в создании эффективных численных методов решения задач выпуклой оптимизации вида ($b = 1, 2$; $l = 1, 2, \infty$)

$$\left\| P^T p - p \right\|_l^b \to \min_{p \geq 0: \sum_{k=1}^{n} p_k = 1}.$$

Эти наработки оказываются полезными, поскольку вектор PageRank можно также понимать как решение такой задачи. Действительно, всегда (из неотрицательности нормы) имеет место неравенство $\|P^T p - p\| \geq 0$, и только на $v$ имеет место равенство $\|P^T v - v\| = 0$.

Из наиболее интересных здесь результатов отметим в частности, оценку метода условного градиента из работы [2] $\|P^T \tilde{p}_T - \tilde{p}_T\|_2 \leq \varepsilon$, которая получается за $C \cdot (n + s^2 \varepsilon^{-2} \ln n)$ арифметических операций. Удивительно в последней оценке то, что в нее не входит $sn$ – число ненулевых элементов матрицы $P$. В частности, при $s \approx \sqrt{n}$ получается сложность $\sim n$, а не $\sim n^{3/2}$, как можно было ожидать (к сожалению, за это есть и плата в виде $\sim \varepsilon^{-2}$). Тем не менее, еще раз повторим, что несмотря на все возможные ускорения вычислений, использование таких методов в наших целях, к сожалению, не представляется возможным.

Отмеченные в этом пункте методы (МПИ и Поляка–Трембы), не исчерпывают линейки методов, гарантирующих малость $\|\tilde{p}_T - v\|$, и в оценку скорости сходимости, которых входит $\alpha$. К таким методам можно отнести, например, метод Д. Спилмана [91], являющийся вариацией описанного выше метода Поляка–Трембы.

Однако мы сосредоточимся далее на, так называемых, методах Монте-Карло [21]. Наряду с тем, что эти методы являются численными методами поиска вектора PageRank, они также позволяют по-новому проинтерпретировать вектор PageRank.

## 4. Markov Chain Monte Carlo и Google problem

В п. 3 мы обратили внимание на то, что формула (1) задает в определенных случаях вполне эффективный численный метод поиска вектора PageRank. Если кратко резюмировать написанное выше, то построенный на основе этой формулы МПИ может быстро сходиться, при этом обеспечивая очень хорошую точность в "правильных категориях", однако каждая его итерация будет стоить $\sim sn$ арифметических операций, что для не очень разреженных постановок задач может оказаться серьезным препятствием. К тому же и задачу PageRank, как правило, надо решать не один раз, а для каждого запроса, искать свой PageRank (детали см. в п. 8), т.е. для эффективной работы реальной поисковой системы остро необходимо не просто за разумное время находить вектор PageRank, но и делать это как можно быстрее, пожертвовав, например, излишней точностью. Жертвовать в МПИ точностью особого смысла не имеет. Ведь, согласно оценке (3), если мы решаем с помощью МПИ задачу поиска вектора PageRank, например, с точностью $\varepsilon \approx 10^{-3}$ за один час, то за два часа работы МПИ мы получим точность $\varepsilon \approx 10^{-6}$. Отсюда видно, что если уж запускать МПИ, то немного подождав, можно получить хорошую точность (плата за это небольшая). Другое дело, если мы попробуем изменить сам метод так, чтобы он был более чувствителен к точности (т.е. в этом смысле хуже себя вел по сравнению с МПИ), но зато "грубо" мог бы быстрее находить вектор PageRank.

Основная идея, лежащая в основе подхода *метода Монте-Карло*, которая в данном контексте получила название *Markov Chain Monte Carlo* (MCMC) [58], заключается в практи-ческом использовании эргодической теоремы из п. 2. То есть, грубо говоря, надо просто запустить человека и достаточно долго подождать, собирая статистику посещенных им районов (web-страниц). При оценке эффективности работы такого метода возникают два вопроса: насколько эффективным можно сделать шаг и сколько шагов надо сделать человеку?

Для ответа на первый вопрос нужно уточнить, в какой ситуации мы находимся. Если из каждой вершины (web-страницы) выходит не более $s \ll n$ ребер, то можно об этом не сильно задумываться, просто каждый раз приготавливая соответствующее распределение вероятностей. Это будет $\sim s$ арифметических операций. При этом память после каждой операции освобождается, т.е. нам не нужно ничего дополнительно к матрице $P$ и текущему вектору частот хранить. Именно в такой ситуации ($s \ll n$) мы и будем находиться, изучая конкретный web-графы, порожденные моделью Бакли–Остгуса (см. пп. 6, 8). Проблемы возникают, когда $s$ является достаточно большим. Тогда можно заранее (препроцессинг) подготовить специальным образом память за $\sim sn$ арифметических операций. Например, поставив в соответствие каждой вершине свое разбиение отрезка $[0,1]$ так, что число подотрезков равно числу выходящих ребер, а длины подотрезков пропорциональны вероятностям, стоящим на ребрах. Тогда на каждом шаге достаточно один раз генерировать равномерно распределенную на отрезке $[0,1]$ случайную величину (стандартные генераторы, как правило, именно это, в первую очередь, и предлагают), и нанести ее на соответствующий (текущей вершине) подотрезок, в зависимости от того, в какой из подотрезков она попала, туда и сдвинуть человека. Минус в таком подходе – использование дополнительно памяти для хранения $\sim sn$ чисел ну и, конечно, затратный препроцессинг. Ведь именно от этого мы и уходили, отказавшись от МПИ. К счастью, тут требуется $\sim sn$ арифметических операций сделать всего один раз (в отличие от МПИ), но все равно, это может быть слишком дорого. Оказывается, при некоторых дополнительных предположениях существуют быстрые способы онлайн разыгрывания направления движения (без существенных затрат памяти и препроцессинга), которые приводят к тому, что шаг выполняется за $\sim \log_2 s \le \log_2 n$ арифметических операций. По-сути, все эти "дополнительные предположения" – это компактное описание формул расчета $\{p_{ij}\}_{j=1}^n$. Если этот набор хранится просто как набор каких-то чисел, то вряд ли что-то можно упростить. Но если мы имеем некоторую модель построения этих чисел, то есть имеется эффективная формула расчета $p_{ij} = p(i,j)$, то можно надеяться и на сложность шага $\sim \log_2 s$. В частности, если имеет место равенство между собой отличных от нуля внедиагональных элементов матрицы $P$ в каждой строке, то сложность шага $\sim \log_2 s$ достигается, например, методом Кнута–Яо [21].

Ответ на второй вопрос (сколько шагов нужно сделать человеку?) в общем случае требует еще больше всяких оговорок (чем ответ на первый вопрос), тем не менее, постараемся грубо сформулировать основной результат [19]. *Пусть $T \gg T_0 \stackrel{def}{=} C\alpha^{-1}\ln(n/\varepsilon)$, $\nu(T)$ – вектор частот (случайный вектор) пребывания в различных вершинах блуждающего человечка после* $T$ *шагов. Тогда с вероятностью не меньше* $1-\sigma$ *(здесь и везде в дальнейшем* $\sigma \in (0,1)$ *— произвольный доверительный уровень) имеет место неравенство*[4]

$$\|v(T) - v\|_2 \le C\sqrt{\frac{\ln n + \ln(\sigma^{-1})}{\alpha T}}, \qquad (7)$$

*где* $\|x\|_2 = \sqrt{\sum_{k=1}^{n} x_k^2}$. Грубость здесь в том, что вместо $\alpha$ правильно писать некоторую другую (как правило, близкую) константу, которую довольно долго определять. Однако во многих важных случаях (в частности, в случае, рассматриваемом далее в данной статье) написанное здесь неравенство верно. Также отметим, что в общем случае вместо $\ln n$ стоит писать $\max_{k=1,\ldots,n} \ln(v_k^{-1})$. Здесь мы воспользовались спецификой нашей задачи, поэтому записали $\ln n$ (детали см. в пп. 6, 8). Оценка (7) позволяет оценить эффективность MCMC.

Приведем сравнительные характеристики описанных (упомянутых) к настоящему моменту методов. "Сложность" понимается как количество арифметических операций типа умножения двух чисел, которые достаточно осуществить, чтобы (в случае MCMC: с вероятностью не меньше $1-\sigma$) достичь точности решения $\varepsilon$ по целевой функции[5] (в Таблице 1 – Цель). Напомним, что вектор $\tilde{p}_T$ – то, что выдает рассматриваемый метод, а $v$ – вектор PageRank (см. (2)). Стоит также пояснить обозначение $\|x\|_\infty = \max_{k=1,\ldots,n} |x_k|$.

Таблица 1. Сравнение свойств методов решения задачи поиска вектора PageRank

| Метод | Сложность | Цель |
|---|---|---|
| МПИ [72] | $\dfrac{sn}{\alpha}\ln\left(\dfrac{2}{\varepsilon}\right)$ | $\|\tilde{p}_T - v\|_1$ |
| Поляка–Трембы [38] | $\dfrac{2sn}{\varepsilon}$ | $\|P^T \tilde{p}_T - \tilde{p}_T\|_1$ |
| Д. Спилмана [91] | $C\cdot\left(n + \dfrac{s^2}{\alpha\varepsilon}\ln\left(\dfrac{1}{\varepsilon}\right)\right)$ | $\|\tilde{p}_T - v\|_\infty$ |

---

[4] Если бы мы имели дело не с марковским процессом, а с совокупностью независимых одинаково распределенных случайных величин, то в правой части оценки (7) в знаменателе под корнем не нужно было бы писать $\alpha$ (см. формулу (9) ниже), что можно понимать как "прореживание" выборки (см. п. 5), с формированием новой (разреженной) выборки с помощью отбора элементов через каждые $\sim \alpha^{-1}$ шагов. Отсюда появляется следующая гипотеза (действительно, имеющая место): после $\sim \alpha^{-1}$ шагов положение блуждающего человека практически не зависит (в вероятностном плане) от стартового положения.

[5] Заметим, что для метода MCMC имеет смысл рассматривать только $\varepsilon \ll n^{-1/2}$, поскольку $\|(n^{-1} \ldots n^{-1})\|_2 = n^{-1/2}$.

| МСМС [19] | $C \cdot \left( n + \dfrac{\log_2 n \cdot \ln(n/\sigma)}{\alpha \varepsilon^2} \right)$ | $\|\tilde{p}_T - v\|_2$ |
|---|---|---|
| вариация метода условного градиента [2] | $C \cdot \left( n + \dfrac{s^2 \ln n}{\varepsilon^2} \right)$ | $\|P^T \tilde{p}_T - \tilde{p}_T\|_2$ |

Сравнение затрудняется тем, что у всех методов разные целевые критерии. Другими словами, у методов разные цели. Также сравнение не может быть полным, если не принимать в расчет то, как параллелится тот или иной метод. Ведь для задач огромных размеров, к коим, безусловно, относится и Google problem, этот вопрос выходит на передний план. Если умножение матрицы на вектор в МПИ хорошо параллелится, то, например, МСМС в описанном здесь варианте не понятно как можно было бы распараллелить.

Решению данной проблемы посвящен следующий пункт, в котором попутно предлагается новый способ интерпретации вектора PageRank, отличный от предложенного С. Брином и Л. Пейджем.

## 5. Параллелизуемый вариант метода МСМС для Google problem. Концепция равновесия макросистемы и теорема Фишера о свойствах оценки максимального правдоподобия

Идею распараллеливания, собственно также как и идеи остальных описанных методов (МПИ, МСМС), подсказывает "жизнь". А именно, в реальном интернете "блуждает" не один человек (пользователь), а много пользователей. Обозначим число независимо блуждающих пользователей через $N$. Рассмотрим сначала (для простоты) ситуацию, когда всего две вершины, т.е. $n = 2$. Поскольку пользователи блуждают независимо, то для каждого из них можно ожидать, что после $T_0 = C\alpha^{-1}\ln(2/\varepsilon)$ шагов вероятность найти его (на шаге $T_0$) в вершине 1 с хорошей точностью ($\sim \varepsilon$) равна $v_1$ (соответственно, в состоянии 2: $v_2 = 1 - v_1$). Это довольно стандартный результат, тесно связанный с оценкой (3). Таким образом, посмотрев на то, в какой вершине находился каждый человек на шаге $T_0$, мы с хорошим приближением получим, так называемую, *простую выборку X* (независимые одинаково распределенные случайные величины) из *распределения Бернулли* с вероятностью успеха (выпадения орла) равной $v_1$. Последнее означает, что каждый человек (подобно монетке) на шаге $T_0$ независимо от всех остальных с вероятностью $v_1$ будет обнаружен в состоянии 1 ("выпал орлом"), а с вероятностью $1 - v_1$ будет обнаружен в состоянии 2 ("выпал решкой"). Чтобы оценить вектор PageRank в данном простом случае достаточно оценить по этой выборке $v_1$ – частоту выпадения орла. Интуиция подсказывает, что в качестве оценки неизвестного параметра $v_1$ следует ис-пользовать $r/N$ – долю людей, которые оказались в состоянии 1. Интуиция подсказывает правильно! Далее будет показано, что такой способ, действительно, оптимальный.

Вероятность, что число людей $r$, которые наблюдались в состоянии 1, равно $k$ может быть посчитана по формуле (*биномиальное распределение*)

$$P(r=k) = C_N^k v_1^k (1-v_1)^{N-k}.$$

Используя (грубый вариант) *формулы Стирлинга* $k! \simeq (k/e)^k$, получаем *оценку типа Санова* [43]

$$C_N^k v_1^k (1-v_1)^{N-k} \simeq \exp(-N \cdot KL(k/N, v_1)).$$

Отсюда по *неравенству Пинскера* [52]

$$KL(p,q) \geq 2(p-q)^2,$$

где

$$KL(p,q) = -p\ln\left(\frac{p}{q}\right) - (1-p)\ln\left(\frac{1-p}{1-q}\right),$$

следует, что *с вероятностью не меньше* $1-\sigma$ имеет место неравенство

$$\left|\frac{r}{N} - v_1\right| \leq C\sqrt{\frac{\ln(\sigma^{-1})}{N}}, \qquad (8)$$

которое иллюстрируется на рис. 4.

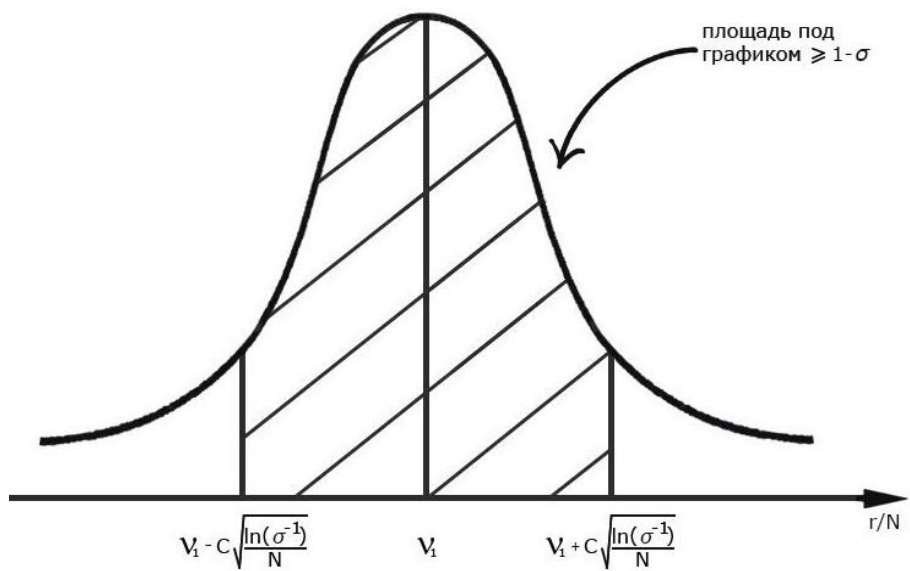

Рис. 4 График зависимости $P(r) = C_N^r v_1^r (1-v_1)^{N-r}$ при большом значении $N$

По-сути, этот рисунок отражает тот факт, что биномиальное распределение (биномиальная мера) с ростом числа исходов (людей) $N$ **концентрируется** вокруг $v_1$.

Вектор $v = (v_1, v_2)^T$, в малой окрестности которого с большой вероятностью на больших временах находится реальный вектор распределения людей по web-страницам (в данном случае, двум), естественно называть **равновесием макросистемы**, описываемой блуждающими по web-графу людьми. *Таким образом, мы пришли к еще одному пониманию вектора Page-Rank $v$, как равновесию описанной выше макросистемы.*

Достаточно большое количество примеров макросистем с анализом их равновесий собрано в главе 6 [10] (см. также [4, 8, 9, 11, 14, 15, 17, 18, 23, 28, 29, 69, 89]). Одним из таких примеров, является "модель Эренфестов" [28, 29], поясняющая различные парадоксы в термодинамике, другим примером макросистемы, которая может быть известна читателям является "Кинетика социального неравенства" [8].

На соотношение (8) можно посмотреть и с другой стороны – с точки зрения математической статистики [32]. А именно, вспомним, что $v_1$ нам не известно. В то время как реализацию случайной величины (статистики) $r/N$ мы можем измерить. Соотношение (8) от такого взгляда на себя не перестанет быть верным.

В общем случае ($n$ вершин) можно провести аналогичные рассуждения: с заменой биномиального распределения *мультиномиальным* [4, 10, 15, 19]. Чтобы получить аналог оценки (8), для **концентрации мультиномиальной меры**, можно использовать векторное неравенство Хефдинга [4, 19, 52, 93]. Оказывается, *с вероятностью не меньше $1-\sigma$ имеет место следующее неравенство [19]*

$$\left\| \frac{r}{N} - v \right\|_2 \le C \sqrt{\frac{\ln(\sigma^{-1})}{N}}, \qquad (9)$$

*где вектор* $r = (r_1, ..., r_n)^T$ *описывает распределение людей по web-страницам в момент наблюдения* $T_0$. Заметим, что константу $C$ здесь можно оценить сверху числом 4. Неравенства (8), (9) являются представителям класса **неравенств концентрации меры**, играющего важную роль в современной теории вероятностей и ее приложениях [24, 52, 73, 93].

Если запустить $N \sim \varepsilon^{-2} \ln(\sigma^{-1})$ человек, дав каждому сделать $T_0 \sim \alpha^{-1} \ln(n/\varepsilon)$ шагов (стоимость шага $\sim \log_2 n$), то (вообще говоря, случайный) вектор $r/N$, характеризующий распределение людей по web-страницам на шаге $T_0$, с вероятностью не меньше $1-\sigma$ обладает таким свойством: $\|r/N - v\|_2 \le \varepsilon$. Это следует из оценки (9). Таким образом, мы приходим к оценке сложности (параллельного варианта) метода MCMC

$$C \cdot \left( n + \frac{\log_2 n \cdot \ln(n/\varepsilon) \ln(\sigma^{-1})}{\alpha \varepsilon^2} \right),$$

которая с точностью до логарифмического множителя совпадает с оценкой, приведенной в п. 4 (см. Таблицу 1). Однако отличие данного подхода в том, что можно пустить блуждать людей параллельно. То есть можно распараллелить работу программы на $N$ процессорах. Разумеется, можно распараллелить и на любом меньшем числе процессоров, давая процессору следить сразу за несколькими людьми. Можно еще немного "выиграть", если сначала запустить меньшее число людей, но допускать, что со временем люди случайно производят клонов, которые с момента рождения начинают жить независимой жизнью, так чтобы к моменту $T_0$ наблюдалось (как и раньше) $N$ человек.

Пока мы не пояснили, почему выбранный способ рассуждения оптимален (и в каком смысле "оптимален"). Нельзя ли как-то точнее/лучше оценить $v$? Чтобы ответить на этот вопрос, снова (для большей наглядности) вернемся к случаю $n = 2$, и заметим, что оценка $r/N$ неизвестного параметра $v_1$ может быть получена следующим образом:

$$r/N = \arg \max_{v \in [0,1]} v^r \cdot (1-v)^{N-r}. \tag{10}$$

Действительно, максимум у функций $f_1(v) = v^r \cdot (1-v)^{N-r_1}$ и

$$f_2(v) = \ln \left( v^r \cdot (1-v)^{N-r} \right) = r \ln v + (N-r) \ln (1-v)$$

достигается в одной и той же точке, поэтому ограничимся задачей $f_2(v) \to \max_{v \in [0,1]}$. Из **принципа Ферма** [45] (максимум лежит либо на границе, либо среди точек экстремума, т.е. среди точек, в которых производная обращается в ноль) получаем условие $df_2(v)/dv = 0$, которое в данном случае примет вид:

$$\frac{r}{v} - \frac{N-r}{1-v} = 0,$$

т.е.

$$\frac{r/N}{v} = \frac{1 - r/N}{1-v}.$$

Отсюда и следует формула (10).

*Функция правдоподобия*, стоящая в правой части (10) отражает вероятность конкретного[6] распределения людей по вершинам (web-страницам), для которого число людей в первой вершине равно $v$. То есть, оценка $r/N$ может быть проинтерпретирована как **оценка макси-**

---
[6] Поскольку "конкретного", то не нужно писать биномиальный коэффициент.

*мального правдоподобия*. Указанный выше способ построения оценок является общим способом получения хороших оценок неизвестных параметров (например, оценки метода наименьших квадратов, которые используются в п. 7, являются оценками максимального правдоподобия, в предположении, что шум имеет нормальное распределение). А именно, пусть есть схема эксперимента (параметрическая модель), согласно которой можно посчитать вероятность того, что мы будем наблюдать то, что наблюдаем. Эта вероятность (плотность вероятности) зависит от неизвестного вектора параметров $\theta$. Будем максимизировать эту вероятность по параметрам этого вектора при заданных значениях наблюдений (выборки). Тогда получим зависимость неизвестных параметров от выборки $X$. Именно эту зависимость $\hat{\theta}^N(X)$ в общем случае и называют *оценкой максимального правдоподобия* вектора неизвестных параметров. Оказывается (**теорема Фишера**), *что в случае выполнения довольно общих условий регулярности используемой параметрической модели такая оценка является асимптотически оптимальной* (**теория ле Кама** [25]). Последнее можно понимать так, что для каждой отдельной компоненты $k$ вектора $\theta$ можно написать неравенство аналогичное (8), справедливое с вероятностью не меньшей $1-\sigma$,

$$\left|\hat{\theta}_k^N(X) - \theta_k\right| \leq C_{k,N}(\theta)\sqrt{\frac{\ln(\sigma^{-1})}{N}}$$

с $C_{k,N}(\theta) \to C_k(\theta)$ при $N \to \infty$ (здесь $N$ – объем выборки, т.е. число наблюдений). При этом $C_k(\theta) > 0$ являются равномерно (по $\theta$ и $k$) наименее возможными. Имеется в виду, что если как-то по-другому оценивать $\theta$ (обозначим другой способ оценивания $\tilde{\theta}^N(X)$), то для любого $\theta$ и $k$ с вероятностью не меньшей $1-\sigma$

$$\left|\tilde{\theta}_k^N(X) - \theta_k\right| > \tilde{C}_{k,N}(\theta)\sqrt{\frac{\ln(\sigma^{-1})}{N}},$$

где $\varliminf_{N \to \infty} \tilde{C}_{k,N}(\theta) \geq C_k(\theta)$.

Строго говоря, именно такая зависимость правой части от $\sigma$ имеет место не всегда. В общем случае при зафиксированном $N$, с уменьшением $\sigma$ правая часть может вести себя хуже, однако при не очень маленьких значениях[7] $\sigma$ написанное верно всегда.

К сожалению, приведенная выше весьма вольная формулировка не отражает в полной мере всю мощь теоремы Фишера. Ведь в таком виде открытым остается вопрос об оценках совместного распределения отклонений сразу нескольких компонент оценки максимального правдоподобия от истинных значений. На самом деле можно сформулировать результат (делается это чуть более громоздко) об асимптотической оптимальности оценки максимального правдоподобия и в таких (общих) категориях. Мы не будем здесь этого делать. Детали см. в классической монографии [25].

---

[7] Пороговое значение $\sigma_0$ удовлетворяет следующей оценке $\ln N \ll \ln \sigma_0^{-1} \ll N$.

Несмотря на, действительно, большую важность сформулированного результата (без преувеличения, являющегося главным результатом математической статистики), стоит сделать несколько оговорок, которые заметно уменьшают весь этот пафос. Главная "проблема" теоремы Фишера (в приведенном нами варианте ле Кама) заключается в том, что оценка оптимальна только в пределе $N \to \infty$. Однако, мы всегда "живем" в условиях конечных выборок (пусть и большого объема). Сформулированный выше результат ничего не говорит, насколько хорошей будет такая оценка при конечном $N$. Также теорема не говорит о том, как получать точные оценки на $C_{k,N}(\theta)$. Теорема лишь предлагает эффективный способ расчета $C_k(\theta)$.[8] Скажем, приведенная нами ранее оценка (8), хотя по форме и выглядит так, как нужно, но все же она далеко не оптимальна. В частности, в оптимальный вариант оценки (8) нужно прописывать под корнем в числителе еще множитель $\nu_1 \cdot (1-\nu_1)$, оцененный нами сверху в (8) константой $1/4$. Мы сделали такую замену в (8), чтобы в правую часть неравенства не входил не известный параметр $\nu_1$. При значениях $\nu_1$ близких к нулю или единице, такая "замена" оказывается слишком грубой. Мы привели простой пример, в котором мы смогли проконтролировать грубость своих рассуждений. В общем случае, к сожалению, это не так просто сделать (если говорить о потерях в мультипликативных константах типа $C$). Поэтому, если мы хотим использовать оценки типа (8), (9), то борьба за "оптимальные константы" сводится не только к выбору оптимального способа оценивания неизвестных параметров, но и к способу оценивания концентрации вероятностного распределения выбранной оценки вокруг истинного значения. За последние несколько лет теория ле Кама была существенно пересмотрена как раз в контексте отмеченных выше проблем. Современный (неасимптотический) вариант параметрической математической статистики, по-прежнему, базирующийся на теореме Фишера, недавно был построен В.Г. Спокойным [92]. Многие отмеченные проблемы удалось в большой степени решить.

Другая проблема (мы ее по ходу уже коснулись) – это зависимость $C_k$ от $\theta$. Ведь нам не известен вектор $\theta$, иначе, зачем тогда его оценивать? Грубо (но практически эффективно) проблема решается подстановкой $C_k(\hat{\theta}^N(X))$. Более точно надо писать неравенство концентрации для $C_k(\hat{\theta}^N(X))$, исходя из неравенства на $\hat{\theta}^N(X)$. Казалось бы, что возникает "порочный круг", и это приводит к иерархии "зацепляющихся неравенств". Однако если выполнить описанное выше огрубление[9] (для того, чтобы обрезать эту цепочку) не сразу, а на каком-то далеком вложенном неравенстве, то чем оно дальше, тем к меньшей грубости это в итоге приведет. Детали мы также вынуждены здесь опустить.

В заключение этого пункта, для закрепления материала, мы предлагаем читателям оценить, сколько надо опросить человек на выходе с избирательных участков большого города, чтобы с точностью 5% (0.05) и с вероятностью 0.99 оценить победителя во втором туре выбо-ров (два кандидата, графы против всех нет). Для решения этой задачи рекомендуется воспользоваться неравенством (8) с явно выписанными константами

$$\left|\frac{r}{N} - v_1\right| \le \frac{1}{2}\sqrt{\frac{\ln(2/\sigma)}{N}}.$$

---
[8] $C_k(\theta)$ – не универсальны, и зависят от использующейся параметрической модели.

[9] Речь идет об оценке $\nu_1 \cdot (1-\nu_1) \leq 1/4$.

## 6. Модель Бакли–Остгуса роста сети Интернет и степенные законы

До сих пор мы считали, что нам дан граф. Однако для ряда приложений (разработка алгоритмов борьбы со спамом и т.п.) необходимо иметь "модель" web-графа, с помощью которой можно было бы исследовать различные закономерности, присущие web-графу. Поясним сказанное. Интернету и многим социальным сетям присущи определенные хорошо изученные закономерности: наличие гигантской компоненты, правило пяти рукопожатий, степенной закон для распределения степеней вершин, специальные свойства кластерных коэффициентов и т.д. [42, 68]. Хотелось бы предложить такую (случайную) модель роста (формирования) этих сетей, которая бы объясняла все эти закономерности. Решив эту задачу, мы можем открывать новые статистические законы, присущие изучаемой сети, с помощью фундаментальной науки, исследуя лишь свойства выбранной модели. Такие исследования позволяют в дальнейшем использовать полученные результаты при разработке алгоритмов уже для реального Интернета.

Одной из лучших на текущий момент моделей Интернета считается модель Бакли–Остгуса [42]. Именно ее мы и будем использовать. А новым законом, который мы хотим изучить, будет степенной закон распределения компонент вектора PageRank, почитанного по графу, сгенерированному по этой модели. Далее мы описываем модель Бакли–Остгуса и выводим[10] степенной закон распределения степеней вершин графа, построенного по этой модели. Именно этот закон и "порождает" степенной закон распределения компонент вектора PageRank.

Более подробно о том, почему "в жизни" так часто возникают степенные законы можно прочитать в обзорах [39, 77, 84, 85].

Итак, рассмотрим следующую модель роста сети.

*База индукции.* Сначала имеется всего одна вершина, которая ссылается сама на себя (вершина с петлей).

*Шаг Индукции.* Предположим, что уже имеется некоторый (ориентированный граф). Пусть появляется новая вершина. Тогда с вероятностью $\beta > 0$ из этой вершины проводится ребро равновероятно в одну из существующих вершин, а с вероятностью $1 - \beta$ из этой вершины проводится ребро в одну из существующих вершин не равновероятно, а с вероятностями пропорциональными входящим степеням вершин[11] – *правило предпочтительного присоединения* (от англ. preferential attachment).

---

[10] Не строго, а в, так называемом, *термодинамическом пределе*. Отметим также, что описываемая модель немного отличается от "настоящей" модели Бакли–Остгуса [42]. В настоящей модели, новая вершина может сослаться и сама на себя. Однако это никак не отразится на основных статистических свойствах построенного по модели случайного графа.

Другими словами, если уже построен граф из $n-1$ вершины, то новая $n$-я вершина сошлется на вершину $i = 1, ..., n-1$ с вероятностью

$$\frac{\operatorname{in deg}_{n-1}(i) + a}{(n-1)(a+1)},$$

где $\operatorname{in deg}_{n-1}(i)$ – входящая степень вершины $i$ в графе, построенном на шаге $n-1$. Параметры $\beta$ и $a$ связаны следующим образом

$$a = \frac{\beta}{1-\beta}.$$

При $a = 1$ получается известная модель Боллобаша–Риордана [42]. Интернету (host-графу) наилучшим образом соответствует значение $a = 0.277$ [42].

Далее вводится число $m$ – среднее число web-страниц на одном сайте, и каждая группа web-страниц с номерами $km+1, km+1, ..., (k+1)m$ объединяется в один сайт. При этом все ссылки (имеющиеся между web-страницами) "наследуются" содержащими их сайтами – получается, что с одного сайта на другой[12] может быть несколько ссылок. Пусть, скажем, получилось, что для заданной пары сайтов таких (одинаковых) ссылок оказалось $l \le m$, тогда мы превращаем их в одну ссылку, но с весом (вероятностью перехода) $l/m$. Именно для так построенного (взвешенного ориентированного) графа мы будем изучать в п. 7 закон распределения компонент вектора PageRank, определяемого формулой (5).

К сожалению, строго доказать, что имеет место степенной закон распределения компонент вектора PageRank по этой модели, насколько нам известно, пока никому не удалось. Имеется только один специальный результат на эту тему, касающийся модели близкой к модели Боллобаша–Риорднана [46]. В данной работе мы также ограничимся только численными экспериментами. Однако, чтобы у читателей появилась некоторая интуиция, почему такой закон может иметь место в данном случае, мы приведем далее некоторые аргументы.

Сначала, следуя работе [77], установим степенной закон распределения входящих вершин в модели Бакли–Остгуса. При этом ограничимся случаем $m = 1$. Обозначим через $X_k(t)$ число вершин с входящей степенью $k$ в момент времени $t$, т.е. когда в графе имеется всего $t$ вершин. Заметим, что по определению

$$t = \sum_{k \ge 0} X_k(t) = \sum_{k \ge 1} k X_k(t) = \sum_{k \ge 0} k X_k(t).$$

---
[11] Выходящая степень всех вершин одинакова и равна 1.
[12] Впрочем, сайты могут совпадать – внутри одного сайта web-страницы также могут друг на друга сослаться.

Поэтому, для $k \geq 1$ вероятность того, что $X_k(t)$ увеличится на единицу при переходе на следующий шаг $t \to t+1$ по формуле полной вероятности равна

$$\beta \frac{X_{k-1}(t)}{t} + (1-\beta)\frac{(k-1)X_{k-1}(t)}{t}.$$

Аналогично, для $k \geq 1$ вероятность того, что $X_k(t)$ уменьшится на единицу при переходе на следующий шаг $t \to t+1$

$$\beta \frac{X_k(t)}{t} + (1-\beta)\frac{kX_k(t)}{t}.$$

Таким образом "ожидаемое" приращение

$$\Delta X_k(t) = X_k(t+1) - X_k(t) \text{ за } \Delta t = (t+1) - t = 1$$

будет[13]

$$\frac{\Delta X_k(t)}{\Delta t} = \beta \frac{(X_{k-1}(t) - X_k(t))}{t} + (1-\beta)\frac{(k-1)X_{k-1}(t) - kX_k(t)}{t}. \qquad (11)$$

Для $X_0(t)$ уравнение аналогичное (11) будет иметь вид

$$\frac{\Delta X_0(t)}{\Delta t} = 1 - \beta \frac{X_0(t)}{t}. \qquad (12)$$

К сожалению, соотношения (11), (12) – не есть точные уравнения, описывающие то, как меняется $X_k(t)$, хотя бы потому, что изменение $X_k(t)$ происходит случайно. Динамика же (11), (12) полностью детерминированная. Однако, для больших значений $t$, когда наблюдается концентрация случайных величин $X_k(t)$ вокруг своих математических ожиданий (средних значений), реальная динамика поведения $X_k(t)$ и динамика поведения средних значений $X_k(t)$ становятся близкими[14] – вариация на тему *теоремы Куртца* [60]. Таким образом, на систему (11), (12) можно смотреть, как на динамику средних значений, вокруг которых плотно сконцентрированы реальные значения. Под плотной концентрацией имеется в виду, что

---

[13] Корректная запись

$$E_{X_{k+1}(t)}\left[\frac{\Delta X_k(t)}{\Delta t} \bigg| X_0(t), ..., X_k(t)\right] = \beta \frac{(X_{k-1}(t) - X_k(t))}{t} + (1-\beta)\frac{(k-1)X_{k-1}(t) - kX_k(t)}{t}.$$

Беря от обеих частей математическое ожидание $E_{X_0(t),...,X_k(t)}[\ ]$, получим

$$E\left[\frac{\Delta X_k(t)}{\Delta t}\right] = \beta \frac{(E[X_{k-1}(t)] - E[X_k(t)])}{t} + (1-\beta)\frac{(k-1)E[X_{k-1}(t)] - kE[X_k(t)]}{t}.$$

[14] Последняя динамика уже является детерминированной динамикой.

разброс значений величины контролируется квадратным корнем из ее среднего значения (см. п. 5).

Будем искать решение системы (11), (12) на больших временах ($t \to \infty$) в виде $X_k(t) \sim c_k \cdot t$ (иногда такого вида режимы называют *промежуточными асимптотиками* [5]). Подставляя это выражение в формулы (11), (12), получим

$$c_0 = \frac{1}{1+\beta}, \quad \frac{c_k}{c_{k-1}} = 1 - \frac{2-\beta}{1+\beta+k\cdot(1-\beta)} \simeq 1 - \left(\frac{2-\beta}{1-\beta}\right)\frac{1}{k}.$$

Откуда получаем следующий **степенной закон**

$$c_k \sim k^{-\frac{2-\beta}{1-\beta}} = k^{-2-a}. \qquad (13)$$

Заметим, что если построить на основе (13) *ранговый закон* распределения вершин по их входящим степеням, т.е. отранжировать вершины по входящей степени, начиная с вершины с самой высокой входящей степенью, то также получим степенной закон [39]

$$\text{in deg}_r \sim r^{-1-\beta}. \qquad (14)$$

Действительно, обозначив для краткости $\text{in deg}_r$ через $x$, получим, что нам нужно найти зависимость $x(r)$, если из формулы (13) известно, что

$$\frac{dr(x)}{dx} \sim -x^{-\frac{2-\beta}{1-\beta}} = -x^{-1-\frac{1}{1-\beta}},$$

где зависимость $r(x)$ получается из зависимости $x(r)$ как решение уравнения $x(r) = x$. Остается только подставить сюда и проверить приведенное соотношение (14). Заметим, что именно ранговый закон распределения компонент вектора PageRank (удовлетворяющего (5) с $\delta > 0$) мы будем численно проверять в п. 7.

Перейдем теперь к пояснению того, почему может иметь место степенной закон распределения компонент вектора PageRank. Для этого предположим, что матрица $P$ имеет вид

$$P \sim \begin{bmatrix} 1^{-\lambda} & 2^{-\lambda} & 3^{-\lambda} & 4^{-\lambda} & 5^{-\lambda} & \dots \\ 1^{-\lambda} & 2^{-\lambda} & 3^{-\lambda} & 4^{-\lambda} & 5^{-\lambda} & \dots \\ 1^{-\lambda} & 2^{-\lambda} & 3^{-\lambda} & 4^{-\lambda} & 5^{-\lambda} & \dots \\ \multicolumn{6}{c}{\dots\dots\dots\dots\dots\dots\dots\dots\dots} \end{bmatrix}.$$

Такой вид матрицы означает, что для каждого сайта имеет место точный (не вероятностный) степенной закон распределения выходящих степеней вершин, имеющий одинаковый вид для всех сайтов. Конечно, это намного более сильное предположение, чем то, что мы выше получили для модели Бакли–Остгуса. Тогда для выписанной матрицы $P$ выполняются условия единственности вектора PageRank $v$, определяющегося формулой (2) (см. п. 2). Более того, этот вектор, неизбежно должен совпадать со строчкой (не важно какой именно – они одинаковые) матрицы $P$, т.е.

$$v_k \sim k^{-\lambda}.$$

Но это и означает, что имеет место степенной закон распределения компонент вектора PageRank. Разумеется, проведенные рассуждения ни в какой степени нельзя считать доказательством. Тем не менее, мы надеемся, что некоторую интуицию эти рассуждения читателям все-таки смогли дать.

В контексте написанного выше хотелось бы отметить, что подобно системе (11), (12) можно записать **динамику средних** (говорят также *квазисредних* [11, 17]) и для макросистемы блуждающих по web-графу людей из предыдущего пункта. А именно, предположим, что людей достаточно много и что каждый человек в (любом) промежутке времени $[t, t+\Delta t)$ независимо от остальных с вероятностью по порядку равной $\Delta t$ совершает переход по одной из случайно выбранных ссылок (согласно матрице $P$). Обозначив через $c_k(t)$ долю людей, находящихся в момент времени $t$ на web-странице с номером $k$ получим следующую систему (читатели, знакомые с дифференциальными уравнениями, здесь могут перейти к пределу при $\Delta t \to 0$)

$$\frac{\Delta c^T(t)}{\Delta t} = c^T(t)P - c^T(t). \tag{15}$$

Формула (15) подтверждает вывод о том, что вектор PageRank $v$, удовлетворяющий системе (2), действительно, можно понимать как равновесие макросистемы. В самом деле, если существует предел $v = \lim_{t \to \infty} c(t)$, то из (15) следует, что этот предел должен удовлетворять (2). Здесь, в отличие от п. 2, предел всегда существует. Но также как и в п. 2 может зависеть от начального условия. Для того чтобы предел не зависел от начального условия (был единственным) нужно сделать предположение о наличии в графе "Красной площади".

Имеется глубокая связь между приведенной выше схемой рассуждений и общими моделями макросистем, которые с точки зрения математики можно понимать как разнообразные модели стохастической химической кинетики [4, 10, 11, 14, 15, 17, 27, 33, 89]. В частности, система (15) соответствует *закону действующих масс Гульдберга–Вааге* [6, 17, 27, 33, 40]. При этом важно подчеркнуть, что возможность осуществлять описанный выше (канонический) *скейлинг*, по-сути заключающийся в замене концентраций (пропорций/долей) их средними значениями, обоснована теоремой Куртца [33, 60] (в том числе и для нелинейных систем – появляющихся, когда имеются не только унарные реакции, как в примере с PageRank'ом) только на конечных отрезках времени. Для бесконечного отрезка[15] требуются дополнительные оговорки, например, выполнение условия детального баланса и его обобщения [4, 6, 10, 11, 17, 33]. Хорошим примером тут является известная модель "хищник–жертва" [9], приводящая к системе Лотки–Вольтерра [41] лишь на конечных отрезках времени. На бесконечном отрезке времени либо сначала все "зайцы" будут съедены "волками", после чего все "волки" погибнут от голода, либо сначала все "волки" погибнут из-за нехватки пищи ("зайцев"), после чего "зайцы" не ограничено расплодятся [10, 17, 23]. И в том и в другом случае такая асимптотика никак не соответствует (нелинейным) незатухающим колебаниям численностей "волков" и "зайцев", которые предписывает решение системы Лотки–Вольтерра. Другими словами, в общем случае использованные нами предельные переходы (по времени и числу агентов) не перестановочны! Рассуждения п. 5, соответствуют следующему порядку предельных переходов: сначала $t \to \infty$ (выходим на инвариантное меру/стационарное распределение), потом $N \to \infty$ (концентрируемся вокруг наиболее вероятного макросостояния инвариантной меры), а рассуждения этого пункта: сначала $N \to \infty$ (переходим на описание макросистемы на языке концентраций, устраняя случайность с помощью законов больших чисел), потом $t \to \infty$ (исследуем аттрактор полученной при скейлинге детерминированное, т.е. не стохастической, системы). Для примера макросистемы из п. 5 получается один и тоже результат. Более того, если посмотреть на то, как именно концентрируется инвариантная мера для этого примера, то получим, что *концентрация экспоненциальная* [4, 10, 17, 43]

$$P(r=k) = \frac{N!}{k_1! \cdot \ldots \cdot k_n!} v_1^{k_1} \cdot \ldots \cdot v_n^{k_n} \simeq \exp(-N \cdot KL(k/N, v)),$$

где *функция действия* (*функция Санова*) $KL(x,y) = -\sum_{k=1}^{n} x_k \ln(x_k/y_k)$. В других контекстах функцию $KL$ чаще называют *дивергенцией Кульбака–Лейблера* или просто *энтропией*.[16] При этом функция $KL(c(t), v)$, как функция $t$, монотонно убывает с ростом $t$ на траекториях системы (15), т.е. является *функцией Ляпунова*. Оказывается, этот факт[17] имеет место и при намного более общих условиях [4, 10, 17, 33]. Точнее говоря, сам факт о том, что *функция, характеризующая экспоненциальную концентрацию инвариантной меры, будет функцией Ляпунова динамической системы, полученной в результате скейлинга из марковского процесса, породившего исследуемую инвариантную меру*, имеет место всегда, а вот то, что именно такая $KL$-функция будет возникать, соответствует макросистемам, удовлетворяющим обобщенному условию детального баланса (условию Штюкельберга–Батищевой–Пирогова [4, 6, 10, 15, 16, 33]), и только таким макросистемам [17].

---

[15] А именно эта ситуация нам наиболее интересна, поскольку, чтобы выйти на равновесие, как правило, необходимо достаточно долго пожить $t \to \infty$.

В заключение этого пункта заметим, что подобно п. 5 можно получить закон (13) в более точных вероятностных категориях. Хотя это можно сделать вполне элементарными комбинаторными средствами (см., например, [65]), тем не менее, соответствующие выкладки оказываются достаточно громоздкие, поэтому мы не приводим их здесь.

---

[16] Отсюда, по-видимому, и пошло, что "равновесие следует искать из *принципа максимума энтропии*" [4, 6, 10, 17, 33, 69, 71].
[17] Известный из курса термодинамики/статистической физики, как *H-теорема Больцмана*.

# 7. Результаты практических экспериментов по поиску вектора Page-Rank для графа, сгенерированного по модели Бакли–Остгуса

В численных экспериментах использовалась модель Бакли–Остгуса (см. п. 6) с разными значениями параметра $a \geq 0$, в которой $m = 10$ (среднее число web-страниц на одном сайте), а число сайтов (далее в этом пункте будем обозначать это число через $n$) было равно $n = 10^5$. Эксперименты проводились на нескольких компьютерах с разными операционными системами (частота процессора 2–3 ГГц, размер оперативной памяти 8–16 Гб), код был реализован на языке Python [98]. Основные временные затраты были связаны с подготовкой web-графа согласно модели п. 6.

Сначала сравнивались методы из Таблицы 1 п. 4 (и не только). Численные эксперименты вполне определенно продемонстрировали, что для web-графа, построенного по модели Бакли–Остгуса (при указанных выше параметрах) при желаемой точности $\varepsilon \ll 10^{-5}$ (см. столбец "Цель" в табл. 1) МПИ доминирует все остальные подходы к расчету вектора PageRank (по всем разумным критериям). В приведенных далее численных экспериментах выбиралась точность $\varepsilon \simeq 10^{-7}$, т.е. компоненты вектора PageRank $v$ восстанавливались следующим образом: $\|\tilde{p}_T - v\|_1 \leq 10^{-7}$. Время работы МПИ (при $a = 1$) составляло порядка 5 минут. С ростом параметра $a$ это время немного увеличивалось.

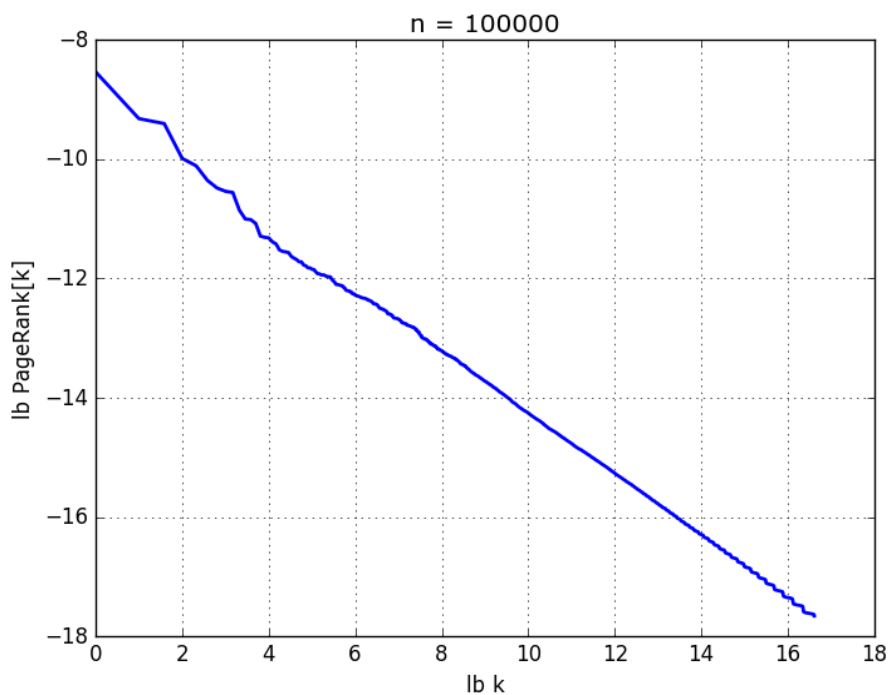

Рис. 5 ($a = 1$) Зависимости отсортированных по убыванию компонент $\log_2 v_k$ от $\log_2 k$

С помощью МПИ была проверена гипотеза, что имеет место степенной закон распределения компонент вектора PageRank, см. рис. 5. При построении графика на рис. 5 компоненты вектора PageRank предварительно сортировались по убыванию, однако, если такую сортиров-ку не производить, то вид графика практически не изменялся. Аналогичные графики получались и на других реализациях web-графа, полученного по модели Бакли–Остгуса, при других значениях параметра $a \geq 0$.

На рис. 6, 7 через $g(a) < 0$ обозначен показатель степени в предполагаемом степенном законе $v_k \sim k^{g(a)}$. Этот показатель определялся по методу наименьших квадратов [67] по данным, отображенным на рис. 5. Для каждого значения $a \geq 0$ по модели Бакли–Остгуса независимо генерировалось 15 web-графов, исходя из разброса посчитанных по этим web-графам значений $g(a)$ около среднего значения, строились приведенные на рис. 6, 7 три кривые (среднее значение/average и крайние значения). Отметим, что при увеличении параметра $m$ в несколько раз $g(a)$ изменялся на несколько тысячных, т.е. можно считать, что $g(a)$ не зависит от выбора $m$. Однако с увеличением $m$ время расчета PageRank'а также увеличивалось. Отметим также, что при $n = 10^4$ получались аналогичные графики только с немного большим разбросом.

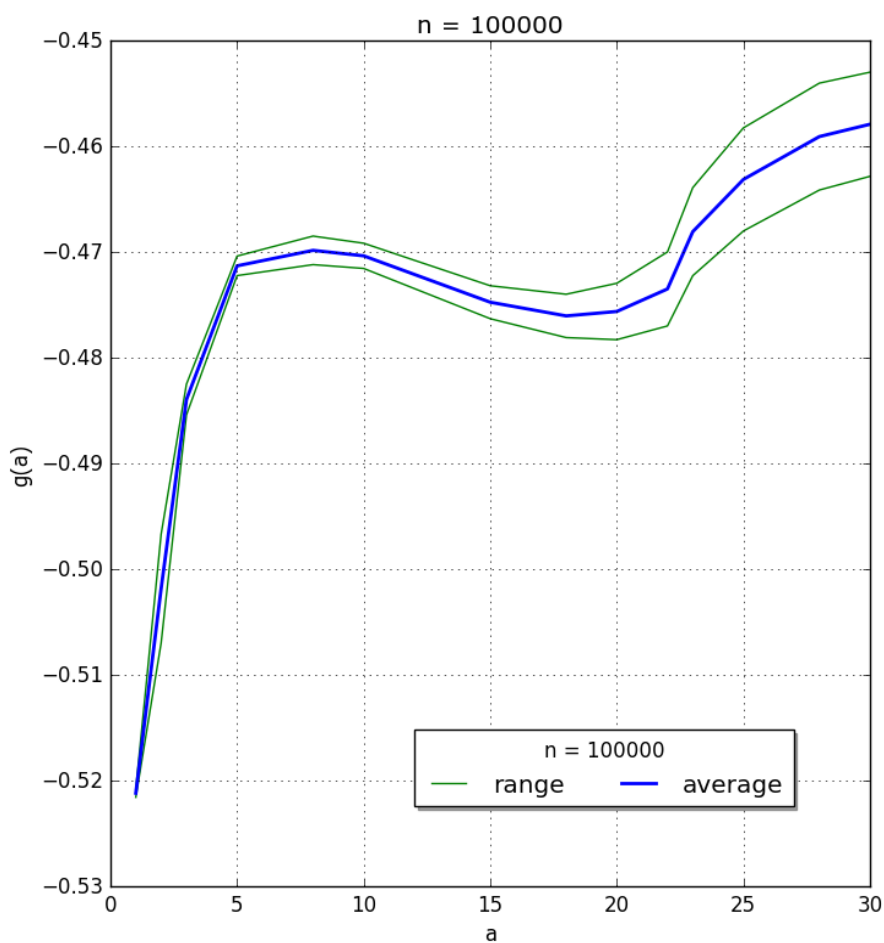

Рис. 6 Зависимость $g(a)$, $a \in [0, 30]$ в предполагаемом законе $v_k \sim k^{g(a)}$

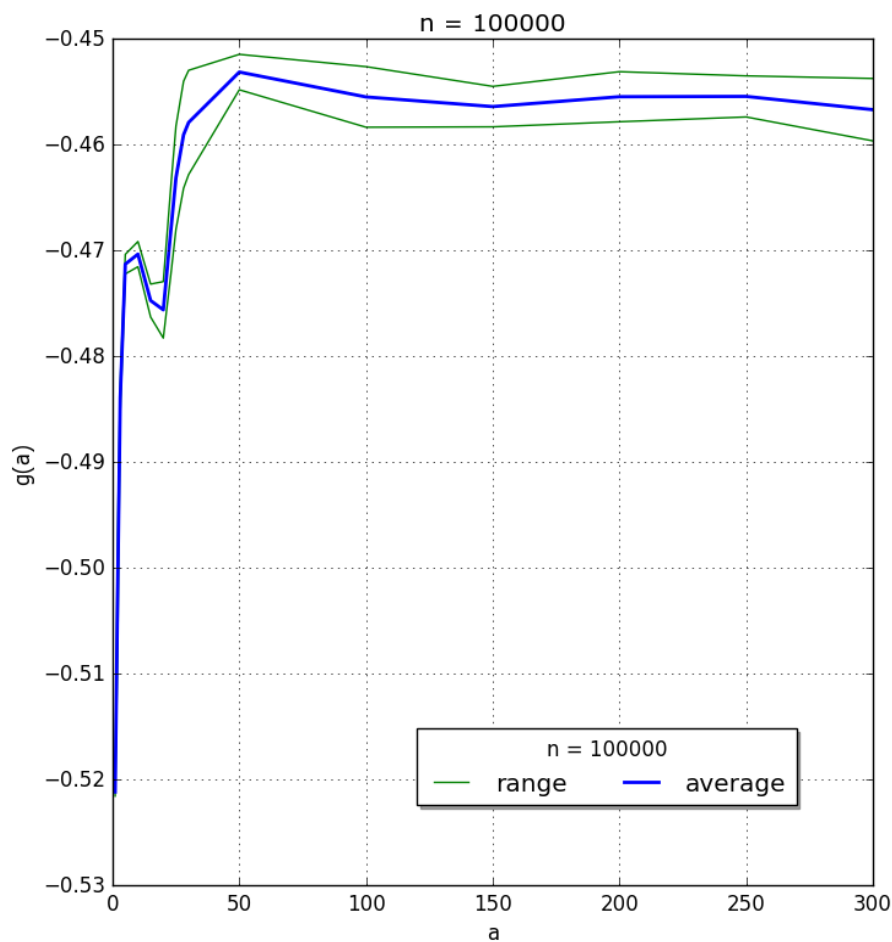

Рис. 7 Зависимость $g(a)$, $a \in [0, 300]$ в предполагаемом законе $\nu_k \sim k^{g(a)}$

## 8. Немного о том, как реально устроены поисковые системы и при чем тут задача поиска вектора PageRank?

Концепция вектора PageRank, описанная в предыдущих пунктах, была нацелена на решение задач ранжирования web-страниц по запросам пользователей. Заметим, что в данной концепции никак не используется информация о заданном запросе. Тем не менее, в этом пункте будет описан один из существующих способов использования именно "классической" концепции вектора PageRank (Брин–Пейдж) для решения реальных задач ранжирования web-страниц.

Может показаться, что в поисковой системе либо есть огромная база соответствий между запросами $q$ и web-страницами $w$, либо есть огромный штат сотрудников, которые за долю секунды успевают просмотреть миллиарды страниц и выбрать нужные. Конечно же, обе гипотезы далеки от истины. Тем не менее, некоторая доля истины есть в обеих. Разумеется, в современных поисковых системах (таких как Яндекс и Google) есть (постоянно обновляемая) огромная база соответствий страниц запросам. Эти соответствия выражаются векторами «признаков» $y(q,w)$. Каждый элемент такого вектора $y_k(q,w)$ (признак) вычисляется как мера обладания пары $(q,w)$ (запрос, web-страница) некоторым, как правило, интуитивным свойством $k$, которым обязана такая пара обладать, если страница должна быть показана в выдаче по запросу (страница *релевантна* запросу). К таким свойствам, например, относят "свежесть", "соответствие тематике" и т.п. Для хорошей работы поисковой системы необходимо много признаков, поэтому векторы признаков в огромных компаниях, занимающихся поиском, имеют, как правило, большую размерность (от 500 и выше), т.е. $k=1,...,l$, $l \sim 10^3$. Кроме того, в компаниях, предлагающих пользователям качественный поиск, работают специально обученные люди, называемые асессорами. Эти специалисты «размечают» (выставляют оценки) некоторые пары (запрос, web-страница) в соответствии со степенью релевантности страницы запросу. Оценки пар выставляются по 5-балльной шкале (от «0» за наименьшую степень релевантности до «4» за наибольшую степень релевантности).

Если задать некоторый запрос поисковой системе, то это еще не значит, что, во-первых, в базе поисковой системы такой запрос имеется, а, во-вторых, что, если такой запрос в базе и имеется, то найдутся web-страницы/документы (или достаточно большое количество документов), размеченные асессорами по запросу. Вообще говоря, вероятность того, что два последних событий выполнены, близка к нулю. В этой связи работники поисковых систем прибегают к решению задачи *машинного обучения* [67, 90], которую (в данном контексте) можно сформулировать следующим образом [51, 54–56].

Прежде всего, опишем модель, положенную в основу ранжирования. Подобно (5') (также с $\delta = 0.15$) для каждого запроса $q$ будем считать, что ранжирование осуществляется согласно вектору $p(q,x)$, являющегося решением системы линейных уравнений (вектор параметров $x$ определяется ниже)

$$p^T = (1-\delta) p^T P(\mathrm{y}(q),x) + \delta \pi^T (\mathrm{y}(q),x), \qquad (16)$$

где $\mathrm{y}(q) = \{y(q,w)\}_{w=1}^n$, а зависимости $P(\mathrm{y}(q),x)$ (матрица) и $\pi(\mathrm{y}(q),x)$ (вектор) считаются известными (каждый элемент этой матрицы или вектора может быть вычислен за $\mathrm{O}(d)$). Если еще знать вектор $x \in \mathbb{R}^d$, $d \simeq 2l \sim 10^3$, то решив (16) МПИ (см. п. 3)

$$p^T(t+1) = (1-\delta) p^T(t) P(\mathrm{y}(q),x) + \delta \pi^T (\mathrm{y}(q),x), \qquad (17)$$

можно осуществить ранжирование. Для определения вектора $x$ используется то, что асессоры для ряда запросов $q \in Q$ разметили ряд документов (web-страниц), т.е., грубо говоря, для запроса $q \in Q$ известны некоторые компоненты "экспертного" вектора PageRank $v^{\exp}(q)$. Считается, что также задана "мера несоответствия" $\mu(p,v)$, согласно которой можно измерить, насколько вектор $p(q,x)$ соответствует вектору $v^{\exp}(q)$. Вектор параметров $x$ определяется исходя из решения следующей задачи оптимизации (обучения)

$$F_Q(x) = \sum_{q \in Q} \mu\big(p(q,x), v^{\exp}(q)\big) \to \min_{x \in \mathbb{R}^d}. \qquad (18)$$

Зависимости $P(\mathrm{y}(q),x)$ и $\pi(\mathrm{y}(q),x)$ можно выбирать с большим произволом. Предположим, что есть два различных сценария 1 и 2. Чтобы определиться с выбором сценария, множество $Q$ разбивается на два подмножества $Q_L$ (обучающее/Learning) и $Q_T$ (контрольное/Test). На обучающем множестве происходит обучение моделей для обоих сценариев: $F_{Q_L}(x) \to \min_{x \in \mathbb{R}^d}$. На выходе получаются векторы $x^1$, $x^2$, отвечающие этим сценариям. Далее в зависимости от того, что меньше $F_{Q_T}(x^1)$ или $F_{Q_T}(x^2)$, отдается предпочтение одному из этих способов ранжирования/сценарию [67]. Для удобства обозначений, в оставшейся части этого пункта будем опускать нижний индекс у $F_Q(x)$.

Зависимости $P(\mathrm{y}(q),x)$, $\pi(\mathrm{y}(q),x)$ и $\mu(p,v)$ стараются выбирать так, чтобы задача (18) была гладкой (добиться еще и выпуклости задачи (18), т.е. выпуклости $F(x)$, к сожалению, пока не удавалось) с липшицевым градиентом:

$$\|\nabla F(y) - \nabla F(x)\|_2 \le L\|y - x\|_2. \qquad (19)$$

В этом случае для решения задачи (18) можно использовать обычный *метод градиентного спуска* (МГС) [35, 86], восходящий к пионерским работам Б.Т. Поляка начала 60-х годов [37]:

$$x^{k+1} = x^k - \frac{1}{L}\nabla F(x^k) \qquad (20)$$

или его *адаптивный вариант*[18] [80]

1. $L^k = L^{k-1}/2$.
2. $x^{k+1} = x^k - \frac{1}{L^k}\nabla F(x^k)$.
3. *Если*[19]

$$F(x^{k+1}) > F(x^k) + \langle \nabla F(x^k), x^{k+1} - x^k \rangle + \frac{L^k}{2}\|x^{k+1} - x^k\|_2^2, \qquad (21)$$

---

[18] Адаптивный вариант МГС является полезным на практике, поскольку значение константы $L$, как правило, неизвестно, а её грубые оценки сверху оказываются завышенными, что замедляет сходимость. Кроме того, в МГС с фиксированным шагом никак не учитывается, что по мере приближения к экстремуму $L$ уменьшается. Если это учитывать, то метод будет быстрее сходиться. В адаптивном МГС происходит настройка на параметр $L$, отвечающий текущему участку пребывания метода (итерационной последовательности), а не на худшую точку, как в обычном МГС. Эксперименты показывают, что за счет адаптивности метод ускоряется на порядок.

[19] Несложно понять, что если в формуле (21) взять $L^k = L$, где $L$ определяется согласно (19), то знак в неравенстве ">" следует поменять на противоположный "≤".

то $L^k := 2L^k$ и возвращение на шаг 2, *иначе* переход на следующую итерацию: $k := k+1$.

Среднее (на одну итерацию) число вычислений значения функции $F(x)$ и градиента в таком методе не превосходит 4 [80].

Используя (19), можно оценить, как сходится (глобально) МГС (20) (аналогично и метод с адаптивным подбором шага) исходя из следующей простой оценки (см., например, [94])

$$F(x^{k+1}) = F\left(x^k - \frac{1}{L}\nabla F(x^k)\right) \leq F(x^k) - \frac{1}{2L}\left\|\nabla F(x^k)\right\|_2^2. \qquad (22)$$

Выберем в качестве *критерия останова метода* условие $\left\|\nabla F(x^k)\right\|_2 \leq \varepsilon$. Тогда, согласно (22), на каждой итерации метода (20) происходит уменьшение целевого функционала не менее чем на $\varepsilon^2/(2L)$. Отсюда можно заключить, что число итераций, которые необходимо сделать до остановки, оценивается как $N \sim 2L/\varepsilon^2$. Несмотря на грубость проведенных рассуждений, оказывается, что оценка

$$N \sim L/\varepsilon^2 \qquad (23)$$

в общем случае является неулучшаемой [34, 63] с точностью до мультипликативной константы, зависящей только от $F(x^0) - F(x_*)$, где $x_*$ – экстремум (т.е. $\nabla F(x_*) = 0$), к которому сходится метод, а $x^0$ – точка старта метода.

Заметим, что МГС (для любой точки старта) сходится к одному из экстремумов функции $F(x)$, вообще говоря, зависящему от точки старта, но, к сожалению, не обязательно к локальному минимуму[20] [35]. Рассчитывать на возможность отыскания глобального минимума $F(x)$ (без дополнительных предположений), к сожалению, не приходится, поскольку для невыпуклых задач даже с единственным локальным минимумом (являющимся глобальным минимумом) существует следующая *нижняя оценка* $N \sim \varepsilon^{-(d-1)}$ на необходимое число итераций [34]. Причем на каждой итерации можно вычислять производные $F(x)$ сколь угодно высокого порядка (если последние существуют) в одной выбранной на этой итерации точке. На практике с отмеченной проблемой часто помогает бороться мультистарт [22]: независимый запуск траекторий метода из разных точек. Однако такой способ поиска глобального минимума является в общем случае очень трудозатратным [22]. Недавние исследования, связанные с *обучением глубоких нейронных сетей*/Deep Learning [66, 90, 96], показали, что совсем не обязательно всегда пытаться найти именно глобальный минимум. Например, если функция (такие функции, по-видимому, часто возникают в Deep Learning и ряде других задач (машинного) обучения) имеет много локальных минимумов приблизительно одной "глубины", то с точки зрения качества обучения не так важно в какой из минимумов в итоге "свалится" метод.

---

[20] Впрочем, путем усложнения описанной процедуры (допуская возможность вычисления гессиана $\nabla^2 F(x^k)$), можно добиться, чтобы сходимость была именно к локальному минимуму, см., например, [94].

Описанный выше подход (МГС и его адаптивный вариант) предполагает возможность вычисления точного градиента $\nabla F(x)$ функции $F(x)$. Однако, итерационный процесс (17) позволяет вычислять (при том только приближенно) лишь значение $p(q,x)$, а следовательно и функции $F(x)$. Если забыть на некоторое время про неточность вычисления $F(x)$, то при весьма общих условиях (см., например, [20, 49, 86]) на основе "графа вычисления" $F(x)$ можно построить "обратный граф" вычисления $\nabla F(x)$ за время не превышающее[21] $4\times$[время расчета $F(x)$], причем в большинстве случаев число 4 может быть уменьшено до числа из интервала $(2,3)$ [49]. Соответствующая общая техника называется[22] (быстрым) *автоматическим дифференцированием* (БАД) [20, 49, 86] (automatic differentiation), а в литературе по нейронным сетям *методом обратного распространения* [49, 96] (back propagation). Интересно при этом заметить, что расчет матрицы Якоби $\left[\partial p(q,x)/\partial x\right]$ не может быть в общем случае осуществлен быстрее, чем за время $2d\times$[время расчета $p(q,x)$].

Вернемся к неточности вычисления $p(q,x)$, а следовательно и $F(x)$. Проблема в том, что в БАД работа должна происходить с точным алгоритмом вычисления значения функции (именно по этом алгоритму строится соответствующая "обратная" процедура расчета градиента). Если алгоритм, вычисляющий функцию $F(x)$, в свою очередь, является неточным, то не понятно, к чему приведет в этом случае техника БАД. Собственно, именно такая проблема и возникает в подходе к решению различных вариационных задач, задач оптимального управления и ряда обратных задач предполагающем сначала дискретизацию задачи (например, за счет замены систему дифференциальных уравнений разностной схемой), а затем использование техники БАД для вычисления градиента дискретизированного функционала [12, 20]. Насколько нам известно, на данный момент в общем случае не существует теоретических обоснований у описанного подхода. Впрочем, в важных частных случаях или в общем случае, но без точных оценок, обоснование имеется в книге [12], см. также цитированную там литературу. Конкретно в рассматриваемом в данном пункте случае за счет специфики процедуры (17) и сделанных предположений удается теоретически обосновать отмеченный переход [51]. А именно, получить оценки на точность вычисления $\nabla F(x)$, исходя из числа сделанных итераций при расчете $p(q,x)$.

Однако можно пойти и по-другому пути, основанному на концепции неточного оракула [16, 37, 51, 57, 59]. Этот путь органично завязывается на то, каким образом можно учитывать влияние неточности в вычислении $F(x)$ и $\nabla F(x)$ на скорость сходимости численного метода для задачи (18) (сейчас будет удобно уже рассматривать сразу адаптивный МГС).

---

[21] Чтобы проще было это понять, можно представить себе, например, функцию $F(x) = \langle c, x \rangle$.
[22] Техника автоматического дифференцирования предполагает гладкость всех функций. В нашем случае это требует некоторых оговорок для существования гладкой зависимости $p(q,x)$ (теорема о неявной функции).

Предположим, что на каждой итерации мы имеем доступ к $(\tilde{\delta}, L)$-оракулу (в нашем случае численной процедуре, построенной на базе (17) – не стоит путать $\delta = 0.15$ в (16), (17) с введенным здесь $\tilde{\delta}$), который для любого $x \in \mathbb{R}^d$ возвращает такие $F_{\tilde{\delta}}(x)$ и $\nabla F_{\tilde{\delta}}(x)$, что

$$\left|F_{\tilde{\delta}}(x) - F(x)\right| \le \tilde{\delta}, \quad \left\|\nabla F_{\tilde{\delta}}(x) - \nabla F(x)\right\|_2 \le \sqrt{8L\tilde{\delta}}$$

и для любого $y \in \mathbb{R}^d$ (несложно показать, что условие "для любого" здесь можно существенно ослабить) имеет место неравенство

$$F(y) \le F_{\tilde{\delta}}(x) + \langle \nabla F_{\tilde{\delta}}(x), y - x \rangle + \frac{L}{2}\|y - x\|_2^2 + \tilde{\delta}.$$

Если адаптивный МГС, с заменой условия (21) на условие

$$F_{\tilde{\delta}}(x^{k+1}) > F_{\tilde{\delta}}(x^k) + \langle \nabla F_{\tilde{\delta}}(x^k), x^{k+1} - x^k \rangle + \frac{L^k}{2}\|x^{k+1} - x^k\|_2^2 + 2\tilde{\delta},$$

и критерием останова $\left\|\nabla F_{\tilde{\delta}}(x^k)\right\|_2 \le \varepsilon/2$, работает с описанным выше $(\tilde{\delta}, L)$-оракулом, где $\tilde{\delta} \simeq \varepsilon^2/(32L)$, то после $N \sim 16L/\varepsilon^2$ итераций $\left\|\nabla F(x^N)\right\|_2 \le \varepsilon$. При доказательстве этой оценки используется то, что всегда $L^k \le 2L$. В общем случае[23] $(\tilde{\delta}, L)$-оракул определяется, как правило, (согласованной) дискретизацией/аппроксимацией процедур расчета $F(x)$ и $\nabla F(x)$. Подчеркнем еще раз, что в предлагаемом подходе $\nabla F_{\tilde{\delta}}(x)$, получается непосредственно из $\nabla F(x)$, а не из $F_{\tilde{\delta}}(x)$. Однако в рассматриваемом в данном пункте примере, как уже отмечалось выше, $\nabla F_{\tilde{\delta}}(x)$ можно получить с нужным теоретическим обоснованием и непосредственно из $F_{\tilde{\delta}}(x)$ с помощью БАД.

С точки зрения практической реализации описанного выше подхода всегда встает вопрос о том, как по заданному $\varepsilon$ строить $(\tilde{\delta}, L)$-оракул, не делая "лишних" вычислений. В рассматриваемом в этом пункте примере такая проблема практически не стоит, поскольку время работы $(\tilde{\delta}, L)$-оракул пропорционально $\sim \ln(\tilde{\delta}^{-1})$ [51]. В частности, в данном случае можно просто пользоваться теоретическими оценками из [51]. К сожалению, в общем случае, возникающие тут теоретические оценки (например, оценки аппроксимации и устойчивости разностной схемы $\tilde{\delta}(h)$, использованной при дискретизации с шагом $h$ системы дифференциальных уравнений в задаче оптимального управления), как правило, оказываются сильно завышенными, поэтому сложность $(\tilde{\delta}, L)$-оракула, пропорциональная $\sim \tilde{\delta}^{-\rho}$, не позволяет надеяться на отсутствие "лишних" вычислений при таком подходе. Разумный выход из данной ситуации заключается в *рестартах* по $h$ [86]: запускаем метод при некотором $h$, затем при $h := h/2$ и т.д. до тех пор, пока при переходе на следующий шаг не будет наблюдаться заметных отличий в результатах. Последнее нуждается в пояснениях (уточнениях), и зависит от специфики рассматриваемой задачи (этому планируется посвятить отдельную работу).

---

[23] Под "общим случаем" понимается всевозможные приложения описанного подхода к решению различных вариационных задач, задач оптимального управления и ряда обратных задач [12, 20, 26]. Отметим, что конечномерность пространства, в котором происходит оптимизация ($d < \infty$), была не существенна (не использовалась при получении оценки на $N$): $x$ может принадлежать бесконечномерному гильбертову пространству. Предлагаемый подход заключается в том, что решать исходную задачу стоит в этом самом пространстве (вообще говоря, бесконечномерном), не дискретизируя исходную постановку, однако при организации вычислительного процесса каждый раз использовать не идеальные значения функции и градиенты, а их приближенные значения. Полученная выше зависимость $\tilde{\delta}(\varepsilon)$ позволяет оптимальным образом подбирать параметры дискретизации/аппроксимации, исходя из желаемой в итоге точности $h(\varepsilon)$.

Выше много внимания было уделено проблеме неточности вычислений, однако, главная проблема БАД в рассматриваемом здесь приложении – это ресурсы памяти. Точнее говоря, проблема большой памяти – это общая проблема БАД [86]: требуется хранить весь граф вычислений значения функции в памяти (со всеми промежуточными выкладками), чтобы можно было построить обратный граф и использовать его для расчета градиента. Однако в нашем случае "большая память" получается уж слишком большой. Действительно, используя обозначения п. 3, можно получить следующую, довольно грубую, оценку на требуемую память (для одного запроса $q \in Q$)

$$\sim sn \cdot \left(d + \ln\left(\varepsilon^{-1}\right)/\delta\right) \geq 10^{13} \cdot 10^{3} \simeq 10^{16},$$

т.е. получается порядка $10^4$ терабайт. Получить доступ к таким ресурсам достаточно быстрой памяти (тем более, оперативной) на деле не представляется возможным с большим запасом. В качестве возможного решения отмеченной проблемы можно предложить использовать вместо МГС его покомпонентный вариант ПМГС (см., например, [16, 79] – в выпуклом случае), который описан ниже в упрощенном варианте:

*равновероятно и независимо от предыдущих розыгрышей выбрать $i_k \in [1,...,n]$,*

$$x_{i_k}^{k+1} = x_{i_k}^{k+1} - \frac{1}{L_{i_k}} \frac{\partial F(x^k)}{\partial x_{i_k}}, \ x_i^{k+1} = x_i^k, \ i \neq i_k,$$

где

$$\left|\partial F(x + he_i)/\partial x_i - \partial F(x)/\partial x_i\right| \leq L_i h, \ e_i - i\text{-й орт.}$$

Анализ скорости сходимости в среднем ПМГС (его адаптивного варианта, а также вариантов этих методов, работающих с неточным оракулом) аналогичен приведенным выше рассуждениям. В частности, можно получить (для упрощения изложения в формуле (24) и в последующем описании критерия останова все огрублено) следующую оценку числа итераций

$$N \sim \bar{L}d/\varepsilon^2 \text{, где } \bar{L} \simeq \frac{1}{n}\sum_{i=1}^{n} L_i \leq L. \tag{24}$$

Отметим, что условие $\|\nabla F(x^k)\|_2 \leq \varepsilon$ теперь уже нельзя проверять на каждой итерации, поскольку тогда теряется смысл в использовании на каждой итерации только одной (случайно выбранной) компоненты градиента, поскольку все равно требуется вычисление полного градиента. Предлагается проверять $\|\nabla F(x^k)\|_2 \leq \varepsilon$ через каждые $\sim d$ итераций. Однако этот способ также упирается в отмеченную выше проблему нехватки памяти, поскольку все равно (пускай и не так часто) приходится считать полный градиент $\nabla F(x^k)$. Другой способ – вычислять $\nabla F(x^k)$ "не честно", используя "модель" $\nabla F(x^k)$, заменяя в которой недоступные на данной $k$-й итерации компоненты градиента последними (на данный момент) известными их значениями.

Несложно построить на базе (17) метод, вычисляющий с требуемой точностью $\partial F(x)/\partial x_i$ за число итераций $\sim \ln(\varepsilon^{-1})/\delta$ (с такой же по порядку стоимостью итерации как у метода (17)) и с затратами памяти $\sim sn \simeq 10$ Тб [51]. Чтобы лучше это понять, можно использовать конечную разность $(F(x+he_i)-F(x))/h$ при специальном выборе $h = h(\varepsilon)$ вместо $\partial F(x)/\partial x_i$ [16, 51, 81]. Все это, конечно, приводит к дополнительным неточностям и, как следствие, замедлению скорости сходимости ПМГС, однако все эти поправки имеют характер логарифмических множителей [51], что не может принципиально изменить качество метода. Таким образом, сопоставляя полученный выигрыш в ресурсах памяти в $d$ раз с потенциально аналогичной потерей в числе итераций (следует сопоставить формулы (23), (24)) можно сказать, что множитель $d$ был просто "перенесен из памяти на время". На самом деле, в (24) может быть $\bar{L} \ll L$. Следуя Ю.Е. Нестерову, рассмотрим следующий пример [16]:

$$F(x) = \frac{1}{2}\langle x, Sx\rangle - \langle b, x\rangle, \tag{25}$$

где $S$ – симметричная матрица, все элементы которой числа от 1 до 2. Тогда

$$L = \lambda_{\max}(S) \geq \lambda_{\max}(1_d 1_d^T) = d \text{, а } L_i = S_{ii} \leq 2 \text{, т.е. } \bar{L} \leq 3.$$

Отметим также, что описанный выше метод (покомпонентный спуск) и его вариации во многих конкретных интересных на практике случаях (см., например, [16]) работают намного эффективнее, чем в рассматриваемом в этом пункте общем случае за счет того, что стоимость пересчета $\partial F(x^{k+1})/\partial x_i$ для ПМГС может быть осуществлена в $\sim d$ раз быстрее, чем расчет $\nabla F(x)$. Это не сложно понять на примере функции (25): полноценный расчет градиента $\nabla F(x)$ стоит $d^2$, а пересчет градиента $\nabla F(x^{k+1}) = Sx^{k+1} - b$ с учетом того, что $Sx^k$ уже известно (с прошлой итерации), и $x_{i_k}^{k+1} = x_{i_k}^{k+1} - L_{i_k}^{-1} \partial F(x^k)/\partial x_{i_k}$, стоит $2d$.

Если в функционале (18) очень много слагаемых ($|Q| \gg 1$) и возможности распараллеливания вычислений ограничены, то использовать описанный выше МГС не разумно из-за огромных вычислительных затрат на каждой итерации. В этом случае у задачи (18) ярко проявляется специальная структура (функционал вида суммы с огромным числом слагаемых), которую можно использовать за счет введения специальных рандомизаций в процедуру МГС [48] (один такой пример рассматривался выше – ПМГС, однако здесь речь идет о специальных рандомизациях суммы).

Для улучшения свойств решения задачи (18) можно вносить различные композиты в функционал [67, 80] (регуляризации, штрафы за сложность модели – попытка контроля переобучения [90]). Описанные выше подходы несложно распространить и на этот случай.

Если задача (рассматриваемая в гильбертовом пространстве)

$$F(x) \to \min_x$$

оказывается выпуклой ($\mu$-сильно выпуклой), то можно говорить о поиске глобального минимума [16, 34, 35, 86]. Будем считать, что на каждой итерации доступен $(\tilde{\delta}, L, \mu)$-оракул [16, 57, 59], который выдает (на запрос, в котором указывается только одна точка $x$) такую пару $(F_{\tilde{\delta}}(x), \nabla F_{\tilde{\delta}}(x))$, что для любых $x, y$ имеет место неравенство

$$\frac{\mu}{2}\|y - x\|^2 \le F(y) - F_{\tilde{\delta}}(x) - \langle \nabla F_{\tilde{\delta}}(x), y - x \rangle \le \frac{L}{2}\|y - x\|^2 + \tilde{\delta}.$$

Существует такая однопараметрическая (параметр $p \in [0,1]$: значение $p = 0$ отвечает обычному методу градиентного спуска [57, 86], а значение $p = 1$ отвечает быстрому градиентному методу Ю.Е. Нестерова [35, 80]) линейка методов (см. [16, 57, 59]), получающих на вход только параметры $L$ и $\mu$, которые работают по следующим (неулучшаемым [57]) оценкам

$$F(x^N) - F_* \le \varepsilon, \ N = \mathrm{O}\left(\min\left\{\left(\frac{LR}{\varepsilon}\right)^{\frac{1}{1+p}}, \left(\frac{L}{\mu}\right)^{\frac{1}{1+p}} \cdot \left\lceil \ln\left(\frac{\mu R^2}{\varepsilon}\right) \right\rceil \right\}\right), \ \tilde{\delta} \le \mathrm{O}\left(\frac{\varepsilon}{N^p}\right), \qquad (26)$$

где можно считать (см. [16]), что $R$ – расстояние от точки старта до решения, а точнее, ближайшего к этой точке решения (если решение не единственно). Если (дополнительно) на одной итерации разрешается обращаться несколько раз (2 – 4) за (приближенным) значением функции $F_{\tilde{\delta}}(x)$, то, подобно написанному ранее в этом пункте, можно предложить адаптивные варианты описываемой линейки методов, которым не требуется подавать на вход параметр $L$ (метод сам настраивается на текущее значение этого параметра), и для которых в оценку (26) входит не худшая (по всем шагам итерационного процесса) константа $L$, а некоторая средняя (пространственно средняя). Если оптимизация происходит в конечномерном пространстве $x \in \mathbb{R}^d$, то при весьма общих условиях [16] можно предложить еще покомпо-нентные варианты описанной адаптивной линейки методов, для которых в оценки числа итераций добавляется множитель $\sim d$, но при этом и стоимость итерации уменьшается также в $\sim d$ раз. Однако при этом константа $L$ усредняется не только по пространству, но и по направлениям всех ортов (не по худшему направлению, как было раньше). Выгода от такой замены также может быть порядка $\sim d$ раз. То есть в целом, если условия позволяют [16], то лучше, в конечном итоге, использовать именно покомпонентные методы. Начиная с работы Ю.Е. Нестерова [79] эти методы стали повсеместно использоваться для решения всевозможных выпуклых задач [16], приходящих из Big Data Science. Особо здесь можно отметить активность P. Richtarik'a [97].

Можно распространить приведенные выше результаты и на случай негладких (не обязательно выпуклых) постановок задач [16, 64, 82]. Также можно исследовать вопрос о том, когда и каким образом стоит рандомизировать [16] описанные выше процедуры, с целью сокращения общего времени работы метода. Интересно также заметить, что теория регуляризация выпуклых постановок задач довольно неплохо разработана к настоящему моменту (см., например, [12, 16, 57]), что позволяет получать сходящиеся по аргументу последовательности (полезное на практике свойство). Выпуклость задачи также обеспечивает возможность замены исходной задачи двойственной к ней. Решение последней задачи в ряде случаев бывает проще осуществить. А с учетом того, что линейка методов (26) – это линейка прямо-двойственных методов [16], то по последовательности, полученной при решении двойственной задачи, можно без существенных дополнительных затрат восстановить с такой же точностью (с которой решалась двойственная задача) и решение исходной (прямой) задачи.

Возвращаясь к затронутой ранее проблеме оптимального выбора шага дискретизации в зависимости от желаемой точности решения исходной задачи $h(\varepsilon)$, возникающего при численном поиске градиента функционала, можно заметить, что для выпуклых постановок задач ситуация заметно интереснее, чем для невыпуклых. А именно, из формулы (26) можно получить зависимость $\tilde{\delta}_p(\varepsilon)$. Пускай имеется оценка $\tilde{\delta}(h)$ (обычно такие оценки как-то можно получать, однако, как уже отмечалось ранее, часто они оказываются завышенными). Наконец, пускай имеется оценка того, сколько стоит итерация $T$, в зависимости от $h$: $T(h)$. Исходя из зависимостей $\tilde{\delta}_p(\varepsilon)$ и $\tilde{\delta}(h)$ можно построить зависимость $h_p(\varepsilon)$. Тогда общее время работы метода будет равно $T(h_p(\varepsilon)) \cdot N_p(\varepsilon)$. Последнее выражение зависит от параметра $p \in [0,1]$. Исходя из минимизации этого выражения по $p \in [0,1]$, можно подобрать оптимальное значение этого параметра. Заметим, что для невыпуклых задач такой степени свободы не было.

В заключение заметим, что в довольно большом числе приложений, в действительности, имеют дело с функционалами вида $F(x) = \frac{1}{2}\|Ax - b\|_2^2$ (см., например, п. 3). В частности, очень популярны такого типа функционалы, при решении обратных задач, в которых $x$ является элементом пространства достаточно гладких функций с заданными краевыми условиями, а $A$ является дифференциальным оператором [26]. Описанные выше довольно общие подходы, как ни странно, хорошо подходят и для решения таких (специальных) задач.